\documentclass[a4paper, 10pt]{amsart}

\usepackage{amsmath,amssymb}

\newtheorem{theorem}{Theorem}[section]
\newtheorem{lemma}[theorem]{Lemma}
\newtheorem{corollary}[theorem]{Corollary}
\newtheorem{proposition}[theorem]{Proposition}

\theoremstyle{definition}

\newtheorem{definition}[theorem]{Definition}
\newtheorem{remark}[theorem]{Remark}

\newtheorem{examples}[theorem]{Examples}

\newcommand{\N}{\mathbb N}
\newcommand{\Z}{\mathbb Z}
\newcommand{\R}{\mathbb R}

\newcommand{\red}{{\text{\rm red}}}

\newcommand{\Fc}{\mathcal F}
\newcommand{\vp}{\mathsf v}

 \DeclareMathOperator{\ord}{ord}

 \DeclareMathOperator{\supp}{supp}

\renewcommand{\t}{\, | \,}
\renewcommand{\time}{\negthinspace \times \negthinspace}

\newcommand{\be}{\begin{equation}}
\newcommand{\ee}{\end{equation}}
\newcommand{\ber}{\begin{eqnarray}}
\newcommand{\eer}{\end{eqnarray}}
\newcommand{\nn}{\nonumber}

\newcommand{\Sum}[2]{\underset{#1}{\overset{#2}{\sum}}}

\numberwithin{equation}{section}

\address{Institut f\"ur Mathematik und Wissenschaftliches Rechnen \\
Karl--Franzens--Universit\"at Graz \\
Heinrich\-stra{\ss}e 36\\
8010 Graz, Austria} \email{alfred.geroldinger@uni-graz.at,
diambri@hotmail.com, wolfgang.schmid@uni-graz.at}

\author{Alfred Geroldinger and David J. Grynkiewicz and Wolfgang A. Schmid}

\begin{document}
\title{The catenary degree  of Krull monoids I}

\thanks{This work was supported by the Austrian Science Fund FWF
(Project Numbers )}

\keywords{non-unique factorizations,  Krull monoids, catenary
degree, zero-sum sequence}

\subjclass[2000]{11R27, 13F05,  20M14}

\begin{abstract}
Let $H$ be a Krull monoid with finite class group $G$ such that
every class contains a prime divisor (for example, a ring of
integers in an algebraic number field or a holomorphy ring in an
algebraic function field). The catenary degree $\mathsf c (H)$ of
$H$ is the smallest integer $N$ with the following property: for
each $a \in H$ and each two factorizations $z, z'$ of $a$, there
exist factorizations $z = z_0, \ldots, z_k = z'$ of $a$ such that,
for each $i \in [1, k]$, $z_i$ arises from $z_{i-1}$ by replacing at
most $N$ atoms from $z_{i-1}$ by at most $N$ new atoms. Under a very
mild condition on the Davenport constant of $G$, we establish a new
and simple characterization of the catenary degree. This
characterization gives a new structural understanding of the
catenary degree. In particular, it clarifies the relationship between $\mathsf c (H)$
and the set of distances of $H$ and opens the way towards obtaining more detailed results on the catenary degree.
As first applications, we give a new upper bound on $\mathsf c(H)$ and characterize when $\mathsf c(H)\leq 4$.
\end{abstract}

\maketitle

\bigskip
\section{Introduction} \label{1}
\bigskip

In this paper we study the arithmetic of Krull monoids, focusing on the case that the class group is finite,
and in addition, we often suppose that every class contains
a prime divisor. This setting includes, in particular, rings of
integers in algebraic number fields and holomorphy rings in
algebraic function fields (more examples are given in Section
\ref{2}). Let $H$ be a Krull monoid with finite class group. Then sets of lengths of
$H$ have a well-defined structure: they are AAMPs (almost
arithmetical multiprogressions) with universal bounds on all
parameters (see \cite[Section 4.7]{Ge-HK06a} for an overview).
Moreover, a recent realization theorem reveals that this description
of the sets of lengths is best possible (see \cite{Sc09a}).

Here we focus on the catenary degree of $H$. This invariant
considers factorizations in a more direct way and not only their
lengths, and thus has found strong attention in the recent
development of factorization theory (see \cite{C-G-L-P-R06,
Ge-Ha08b, C-G-L09, Ge09a, B-C-R-S-S10}). The catenary degree
$\mathsf c (H)$ of $H$ is defined as the smallest integer $N$ with
the following property: for each $a \in H$ and each two
factorizations $z$ and $z'$ of $a$, there exist factorizations $z = z_0,
\ldots, z_k = z'$ of $a$ such that, for each $i \in [1, k]$, $z_i$
arises from $z_{i-1}$ by replacing at most $N$ atoms from $z_{i-1}$
by at most $N$ new atoms. The definition reveals immediately that
$H$ is factorial  if and only if its catenary degree equals zero.
Furthermore, it is easy to verify that the finiteness of the class
group implies the finiteness of the catenary degree, and that the
catenary degree depends only on the class group (under the
assumption that every class contains a prime divisor). However,
apart from this straightforward information, there is up to now
almost no insight into the structure of the concatenating chains of
factorizations and no information on the relationship between the
catenary degree and other invariants such as  the set of distances.
Almost needless to say, apart from very simple cases, the
precise value of the catenary degree---in terms of the group
invariants of the class group---is unknown.

The present paper  brings some light into the nature of the catenary
degree. To do so, we introduce a new arithmetical invariant, $\daleth
(H)$, which is defined as follows (see Definition \ref{3.1}): for
each two atoms $u, v \in H$, we look at a factorization having the
smallest number of factors besides two, say $uv = w_1 \cdot \ldots
\cdot w_s$, where $s \ge 3$, $w_1, \ldots, w_s$ are atoms of $H$ and
$uv$ has no factorization of length $k$ with $2 < k < s$. Then
$\daleth (H)$ denotes the largest possible value of $s$ over all
atoms $u, v \in H$. By definition, we have $\daleth (H) \le \mathsf
c (H)$, and Examples \ref{3.3} offer a list of well-studied monoids
where $\daleth (H)$ is indeed strictly smaller than $\mathsf c (H)$.
But the behavior is different for Krull monoids $H$ with finite class group
and every class containing a prime divisor. Under a very mild
condition on the Davenport constant of the class group, we show that
the catenary degree is equal to $\daleth (H)$ (see Corollary
\ref{4.3} and Remark \ref{4.4}), which immediately implies that the
catenary degree equals the maximum of the set of distances plus two.

Since $\daleth (H)$ is a much more accessible invariant than the
original condition given in the definition of the catenary degree,
the equality $\daleth (H) = \mathsf c (H)$ widely opens the door for
further investigations of the catenary degree, both for explicit
computations as well as for more abstract studies based on methods
from Additive and Combinatorial Number Theory (the latter is done in
\cite{Ge-Gr-Sc11b}, with a focus on groups with large exponent). Exemplifying this, in
Section \ref{5}, we derive an upper bound on $\daleth (H)$, and thus on $\mathsf c(H)$ as well, and then characterize Krull monoids with small catenary degree
(Corollary \ref{5.6}).

\bigskip
\section{Preliminaries} \label{2}
\bigskip

Our notation and terminology are consistent with \cite{Ge-HK06a}. We
briefly gather some key notions. We denote by $\N$ the set of
positive integers, and we put $\N_0 = \N \cup \{0\}$. For real
numbers $a, b \in \R$, we set $[a, b] = \{ x \in \Z \mid a \le x \le
b \}$, and we define $\sup \emptyset = \max \emptyset = \min
\emptyset = 0$. By a \ {\it monoid}, \ we always mean a commutative
semigroup with identity which satisfies the cancellation law (that
is, if $a,b ,c$ are elements of the monoid with $ab = ac$, then $b =
c$ follows). The multiplicative semigroup of non-zero elements of an
integral domain is a monoid.

Let $G$ be an additive abelian group and $G_0 \subset G$ a subset.
Then $[G_0] \subset G$ denotes the submonoid generated by $G_0$ and
$\langle G_0 \rangle \subset G$ denotes the subgroup generated by
$G_0$. We set $G_0^{\bullet} = G_0 \setminus \{0\}$. A family
$(e_i)_{i \in I}$ of \emph{nonzero} elements of $G$ is said to be \
{\it independent} \ if
\[
\sum_{i \in I} m_ie_i =0 \quad \text{implies} \quad m_i e_i =0
\quad \text{for all }  i \in I, \quad \mbox{ where } m_i\in \Z\,.
\]
If $I = [1, r]$ and $(e_1, \ldots, e_r)$ is
independent, then we simply say that $e_1, \ldots, e_r$ are
independent elements of $G$.  The tuple $(e_i)_{i \in I}$ is called a {\it basis} if
$(e_i)_{i \in I}$ is independent and $\langle \{e_i \mid i \in I \}
\rangle = G$.

Let \ $A,\, B \subset G$ \ be subsets. Then \ $A + B = \{ a
+ b \mid a \in A, b \in B \}$ \ is their \ {\it sumset}. If $A
\subset \Z$, then the {\it set of distances} of $A$, denoted $\Delta (A)$, is  the set of all differences between consecutive elements of $A$, formally, all $d \in \N$ for which there exist $l \in
A$ such that $A \cap [l, l+d] = \{l, l+d\}$.  In particular, we have
$\Delta ( \emptyset ) = \emptyset$.

For $n \in \mathbb N$, let $C_n$ denote a cyclic group with $n$
elements. If $G$ is finite with $|G| > 1$, then we have
\[
G \cong C_{n_1} \oplus \ldots \oplus C_{n_r}  \,, \quad \text{and we
set} \quad \mathsf d^* (G) = \sum_{i=1}^r (n_i-1) \,,
\]
where $r = \mathsf r (G) \in \mathbb N$ is the \ {\it rank} of
$G$, $n_1, \ldots , n_r \in \mathbb N$ are integers with  $1 < n_1
\mid \ldots \mid n_r$ and $n_r = \exp (G)$ is the exponent of $G$.
If $|G| = 1$, then $\mathsf r (G) = 0$, $\exp (G) = 1$, and $\mathsf d^* (G) = 0$.

\medskip
\noindent
{\bf Monoids and factorizations.} Let $H$ be a monoid. We
denote by $H^{\times}$ the set of invertible elements of $H$, and we
say that $H$ is \ {\it reduced} \ if $H^{\times} = \{1\}$. Let
$H_{\red} = H/H^{\times} = \{ aH^{\times} \mid a \in H \}$ be the
associated reduced monoid and $\mathsf q (H)$ a quotient group of
$H$. For a subset $H_0 \subset H$, we denote by $[H_0] \subset H$
the submonoid generated by $H_0$. Let \ $a, b \in H$. We say that \
$a$ {\it divides} $b$ \ (and we write $a \t b$) if there is an
element $c \in H$ such that $b = ac$, and we say that \ $a$ and $b$
\ are {\it associated} ($a \simeq b$) if $a \t b$ and $b \t a$.

A monoid $F$ is called \ {\it free $($abelian, with basis $P \subset
F)$} \ if every $a \in F$ has a unique representation of the form
\[
a = \prod_{p \in P} p^{\mathsf v_p(a) } \quad \text{with} \quad
\mathsf v_p(a) \in \N_0 \ \text{ and } \ \mathsf v_p(a) = 0 \ \text{
for almost all } \ p \in P \,.
\]
We set $F = \mathcal F(P)$ and   call
\[
|a|_F = |a| = \sum_{p \in P} \mathsf v_p (a) \quad \text{the \ {\it
length}} \ \text{of} \ a \,.
\]
We denote by $\mathcal A (H)$  the  {\it set of atoms}  of $H$, and
we call $\mathsf Z (H) = \mathcal F ( \mathcal A (H_{\red}))$ the
{\it factorization monoid} of $H$. Further, $\pi \colon
\mathsf{Z}(H) \to H_{\text{red}}$ denotes the natural homomorphism given by mapping a factorization to the element it factorizes.
For $a \in H$, the set
\[
\begin{aligned}
\mathsf Z (a) & = \mathsf Z_H (a)   = \pi^{-1} (aH^\times) \subset
\mathsf Z (H) \quad \text{is called the  {\it set of factorizations}
of
 $a$},
 \\
\mathsf L (a) & = \mathsf L_H (a)   = \bigl\{ |z| \, \bigm| \, z \in
\mathsf Z (a) \bigr\} \subset \N_0 \quad \text{is called the {\it
set of lengths}  of $a$,} \text{ and} \\
\Delta (H) & = \bigcup_{a \in H} \Delta \bigl( \mathsf L (a) \bigr)
\ \subset \N \qquad \text{denotes the {\it set of distances} of $H$}
\,.
\end{aligned}
\]
The monoid $H$ is called
\begin{itemize}
\smallskip
\item {\it atomic} \ if $\mathsf Z (a) \ne \emptyset$ for all $a \in
      H$ (equivalently, every non-unit of $H$ may be written as a finite
      product of atoms of $H$).

\smallskip
\item {\it factorial} \ if $|\mathsf Z (a)| = 1$ for all $a \in
      H$ (equivalently, every non-unit of $H$ may be written as a finite
      product of primes of $H$).
\end{itemize}
Two factorizations  $z,\, z' \in \mathsf Z (H)$  can be written in
the form
\[
z = u_1 \cdot \ldots \cdot u_lv_1 \cdot \ldots \cdot v_m \quad
\text{and} \quad z' = u_1 \cdot \ldots \cdot u_lw_1 \cdot \ldots
\cdot w_n
\]  with
\[
\{v_1 ,\ldots, v_m \} \cap \{w_1, \ldots, w_n \} = \emptyset,
\]
where  $l,\,m,\, n\in \N_0$ and $u_1, \ldots, u_l,\,v_1, \ldots,v_m,\,
w_1, \ldots, w_n \in \mathcal A(H_\red)$. Then $\gcd (z, z') = u_1 \cdot \ldots \cdot u_l$, and we call
$\mathsf d (z, z') = \max \{m,\, n\} = \max \{ |z \gcd (z,
z')^{-1}|, |z'\gcd (z, z')^{-1}| \}\in \N_0$ the {\it distance}
between $z$ and $z'$.

\medskip
\noindent {\bf Krull monoids.} The theory of Krull monoids is
presented in the monographs \cite{HK98, Gr01, Ge-HK06a}. We briefly
summarize what is needed in the sequel. Let $H$ and $D$ be
monoids. A monoid homomorphism $\varphi \colon H \to D$ \ is called
\begin{itemize}
\smallskip
\item a  {\it divisor homomorphism} if $\varphi(a)\mid\varphi(b)$ implies  $a \t b$,  for all $a,b \in H$.

\smallskip
\item  {\it cofinal} \ if, for every $a \in D$, there exists some $u
      \in H$ such that $a \t \varphi(u)$.

\smallskip
\item  a {\it divisor theory} (for $H$) if $D = \mathcal F (P)$
for some set $P$, $\varphi$ is a divisor homomorphism, and for every
$p \in P$ (equivalently, for every $a \in \mathcal{F}(P)$), there
exists a finite subset $\emptyset \ne X \subset H$ satisfying $\gcd \bigl( \varphi(X) \bigr)=p$.
\end{itemize}
Note that, by definition, every divisor theory is cofinal. We call
$\mathcal{C}(\varphi)=\mathsf q (D)/ \mathsf q (\varphi(H))$ the
class group of $\varphi $ and use additive notation for this group.
For \ $a \in \mathsf q(D)$, we denote by \ $[a] = [a]_{\varphi} = a
\,\mathsf q(\varphi(H)) \in \mathsf q (D)/ \mathsf q (\varphi(H))$ \
the class containing \ $a$.  If $\varphi \colon H \to \mathcal F
(P)$ is a cofinal divisor homomorphism, then
\[
G_P = \{[p] = p \mathsf q (\varphi(H)) \mid p \in P \} \subset
\mathcal{C}(\varphi)
\]
is called the \ {\it  set of classes containing prime divisors}, and
we have $[G_P] = \mathcal{C}(\varphi)$. If $H \subset D$ is a
submonoid, then $H$ is called \emph{cofinal} (\emph{saturated},
resp.) in $D$ if the imbedding $H \hookrightarrow D$ is cofinal (a
divisor homomorphism, resp.).

The monoid $H$ is called a {\it Krull monoid} if it satisfies one of
the following equivalent conditions (\cite[Theorem
2.4.8]{Ge-HK06a}){\rm \,:}
\begin{itemize}
\item $H$ is $v$-noetherian and completely integrally closed.

\smallskip
\item $H$ has a divisor theory.

\smallskip
\item $H_{\red}$ is a saturated submonoid of a free monoid.
\end{itemize}
In particular, $H$ is a Krull monoid if and only if $H_{\red}$ is a
Krull monoid. Let $H$ be a Krull monoid. Then a divisor theory
$\varphi \colon H \to \mathcal F (P)$ is unique up to unique
isomorphism. In particular, the class group $\mathcal C ( \varphi)$
defined via a divisor theory of $H$ and the subset of classes
containing prime divisors depend only on $H$. Thus it is called the
{\it class group} of $H$ and is denoted by $\mathcal C (H)$.

An integral domain $R$ is a Krull domain if and only if its
multiplicative monoid $R \setminus \{0\}$ is a Krull monoid, and a
noetherian domain is Krull if and only if it is integrally closed.
Rings of integers, holomorphy rings in algebraic function fields, and
regular congruence monoids in these domains are Krull monoids with finite
class group such that every class contains a prime divisor
(\cite[Section 2.11]{Ge-HK06a}). Monoid domains and power series
domains that are Krull are discussed in \cite{Gi84, Ki01, Ki-Pa01}.

\medskip
\noindent {\bf Zero-sum sequences.} Let  \ $G_0 \subset G$ \ be a
subset and $\mathcal F (G_0)$ the free monoid with basis $G_0$.
According to the tradition of combinatorial number theory, the
elements of $\mathcal F(G_0)$ are called \ {\it sequences} over \
$G_0$. For a sequence
\[
S = g_1 \cdot \ldots \cdot g_l = \prod_{g \in G_0} g^{\mathsf v_g
(S)} \in \mathcal F (G_0) \,,
\]
we call $\vp_g(S)$ the {\it multiplicity} of $g$ in $S$,
\[
\begin{aligned}
|S| & = l = \sum_{g \in G} \mathsf v_g (S) \in \mathbb N_0 \
\text{the \ {\it length} \ of \ $S$} \,,  \quad \supp (S) = \{g \in
G \mid \mathsf v_g (S) > 0 \} \subset G \ \text{the \ {\it support}
\ of \
$S$} \,, \\
\sigma (S) & = \sum_{i = 1}^l g_i \ \text{the \ {\it sum} \ of \
$S$} \quad  \text{and} \quad \Sigma (S) = \Big\{ \sum_{i \in I} g_i
\mid \emptyset \ne I \subset [1,l] \Big\} \ \text{ the \ {\it set of
subsums} \ of \ $S$} \,.
\end{aligned}
\]
The sequence $S$ is called
\begin{itemize}
\item {\it zero-sum free} \ if \ $0 \notin \Sigma (S)$,

\item a {\it zero-sum sequence} \ if \ $\sigma (S) = 0$,

\item a {\it minimal zero-sum sequence} \ if it is a nontrivial
      zero-sum
      sequence and every proper  subsequence is
      zero-sum free.
\end{itemize}
The monoid
\[
\mathcal B(G_0) = \{ S \in \mathcal F(G_0) \mid \sigma (S) =0\}
\]
is called the \ {\it monoid of zero-sum sequences} \ over \ $G_0$,
and we have \ $\mathcal B(G_0) = \mathcal B(G) \cap \mathcal
F(G_0)$. Since $\mathcal B (G_0) \subset \mathcal F (G_0)$ is
saturated, $\mathcal B  (G_0)$ is a Krull monoid (the atoms are
precisely the minimal zero-sum sequences). Its significance for the
investigation of general Krull monoids is demonstrated by Lemma
\ref{3.6}.

\smallskip
For every arithmetical invariant \ $*(H)$ \ defined for a monoid
$H$, we write $*(G_0)$ instead of $*(\mathcal B(G_0))$. In
particular, we set \ $\mathcal A (G_0) = \mathcal A (\mathcal B
(G_0))$ and $\Delta (G_0) = \Delta ( \mathcal B (G_0))$. We define
the \ {\it Davenport constant} \ of $G_0$ by
\[
\mathsf  D (G_0) = \sup \bigl\{ |U| \, \bigm| \; U \in \mathcal A
(G_0) \bigr\} \in \N_0 \cup \{\infty\} \,,
\]
and the following properties will be used throughout the manuscript
without further mention. If $G_0$ is finite, then $\mathsf D (G_0)$
is finite (\cite[Theorem 3.4.2]{Ge-HK06a}).  Suppose that $G_0 = G$
is finite. Then
\be\label{soitis}
1 + \mathsf d^* (G) \le  \mathsf D (G)  \,,
\ee
and equality holds if $G$ is a $p$-group or $\mathsf r (G) \le 2$
(see \cite[Chapter 5]{Ge-HK06a} and \cite[Section 4.2]{Ge09a}).

\bigskip
\section{The catenary degree and its refinements} \label{3}
\bigskip

We recall the definition of the catenary degree $\mathsf c (H)$ of
an atomic monoid $H$ and introduce, for all $k \in \N$,  the
refinements $\mathsf c_k (H)$.

\medskip
\begin{definition} \label{3.1}
Let $H$ be an atomic monoid and $a \in H$.
\begin{enumerate}
\item Let $z,\, z' \in \mathsf Z (a)$ be factorizations of $a$ and $N \in
      \N_{\ge 0} \cup \{ \infty\}$. A finite sequence $z_0,\, z_1, \ldots,
      z_k$ in \ $\mathsf Z(a)$ is called an \ $N$-{\it chain of
      factorizations} \ from $z$ to $z'$ \ if $z = z_0$, \ $z' = z_k$ and
      \ $\mathsf d (z_{i - 1}, z_i) \le N$ for every $i \in [1, k]$.

      \smallskip
      \noindent
      If there exists an $N$-chain of factorizations from $z$ to $z'$, we
      say that $z$ and $z'$ \ can be \ {\it concatenated} \ by an
      $N$-chain.

\smallskip
\item Let \ $\mathsf c_H (a) = \mathsf c (a) \in \N _0 \cup
      \{\infty\}$ denote the smallest $N \in \N _0 \cup \{\infty\}$ such
      that any two factorizations $z,\, z' \in \mathsf Z (a)$ can be
      concatenated by an $N$-chain.

\smallskip
\item For $k \in \N$, we set
      \[
      \mathsf c_k (H) = \sup \{ \mathsf c (a) \mid a \in H \ \text{with} \
      \min \mathsf L (a) \le k \} \in \N_0 \cup \{\infty\} \,,
      \]
      and we call
      \[
      \mathsf c (H) = \sup \{ \mathsf c (a) \mid a \in H \} \in \N_0 \cup
      \{\infty\}
      \]
      the \ {\it catenary degree} \ of $H$.

\smallskip
\item We  set
      \[
      \daleth (H) = \sup \big\{\min \bigl( \mathsf L (uv) \setminus \{2\} \bigr) \mid u,\,v \in \mathcal A(H) \big\}
            \,,
      \]
      with the convention that $\min \emptyset = \sup \emptyset =
      0$.
\end{enumerate}
\end{definition}

\medskip
Let all notations be as above. Then $\daleth (H) = 0$ if and only if
$\mathsf L (uv) = \{2\}$ for all $u, v \in \mathcal A (H)$. By
definition, we have $\mathsf c (a) \le \sup \mathsf L (a)$. Let $z,
z' \in \mathsf Z (a)$. Then, by definition of the distance, we have
$z = z'$ if and only if $\mathsf d (z, z') = 0$.  Thus, $\mathsf c
(a) = 0$ if and only if $a$ has unique factorization (that is,
$|\mathsf Z (a)| = 1$), and hence $H$ is factorial if and only if
$\mathsf c (H) = 0$. Suppose that $H$ is not factorial. Then there
is a $b \in H$ having two distinct factorizations $y, y' \in \mathsf
Z (b)$. A simple calculation (see \cite[Lemma 1.6.2]{Ge-HK06a} for
details) shows that
\begin{equation}
2 + \bigl| |y |-|y'| \bigr| \le \mathsf d (y, y') \,, \quad
\text{and hence} \quad  2 + \sup \Delta (\mathsf L(b)) \le \mathsf c
(b)   \,. \label{3.0}
\end{equation}

The following lemma gathers some simple properties of the invariants
introduced in Definition \ref{3.1}.

\medskip
\begin{lemma} \label{3.2}
Let $H$ be an atomic monoid.
\begin{enumerate}
\item We have $0 = \mathsf c_1 (H) \le \mathsf c_2 (H) \le
      \ldots $ and
      \[
      \mathsf c (H) = \sup \{ \mathsf c_k (H) \mid k \in \N \} \,.
      \]

\smallskip
\item We have $\mathsf c (H) = \mathsf c_k (H)$ for all $k \in \N$
      with $k \ge \mathsf c (H)$.

\smallskip
\item If $\mathsf{c}_{k}(H)>\mathsf{c}_{k-1}(H)$ for some $k \in \N_{\ge 2}$, then $\mathsf{c}_{k}(H) \ge k$.

\smallskip
\item $\sup \Delta (H) \le \sup \{ \mathsf c_k (H) - k \mid k \in \N
      \ \text{with} \ 2 \le k < \mathsf c (H) \}$. Moreover, if
      $\mathsf c (H) \in \N$, then there is some minimal $m \in
      \N$ with $\mathsf c (H) = \mathsf c_m (H)$, and then
      \[
      \sup \{ \mathsf c_k (H) - k \mid k \in \N_{\geq 2}
       \} = \max \{ \mathsf
      c_{k} (H) - k \mid k \in [2,m] \} \,.
      \]

\smallskip
\item For every $k \in \N$, we have
      \[
      \begin{aligned}
      \mathsf c_k (H) & \ge  \sup \{ \mathsf c (a) \mid a \in H \
      \text {with} \ k \in \mathsf L (a) \} \\
        & \ge \sup \{  \mathsf c (a) \mid a \in H \
      \text {with} \ k = \min \mathsf L (a) \} \,,
      \end{aligned}
      \]
      and equality holds if $H$ contains a prime element.

\smallskip
\item If $H$ is not factorial, then
      \begin{equation}
      \daleth (H) \le \min \big\{2 + \sup \Delta (H)  \,, \mathsf c_2 (H) \big\}
      \le \max \bigl\{2 + \sup \Delta (H)  \,, \mathsf c_2 (H) \bigr\} \le
      \mathsf c (H) \,. \label{crucial-inequ}
      \end{equation}
\end{enumerate}
\end{lemma}

\begin{proof}
1. Obvious.

\smallskip
2. If $\mathsf c (H)$ is either zero or infinite, then the assertion
is clear. Suppose that $\mathsf c (H) = m \in \N$. Then there is an $a
\in H$ with factorizations $z = u_1 \cdot \ldots \cdot u_l \in \mathsf Z
(a)$ and $z' = v_1 \cdot \ldots \cdot v_m \in \mathsf Z (a)$, where
$l \in [1, m]$ and $u_1, \ldots, u_l, v_1, \ldots, v_m \in \mathcal
A (H_{\red})$, such that $\mathsf d (z, z') = \max \{l, m\} = m$ and
$z$ and $z'$ cannot be concatenated by a $d$-chain of factorizations
for any $d < m$. Since $\min \mathsf L ( a ) \le m$, we get, for all
$k \ge m$, that
\[
m \le \mathsf c (a) \le \mathsf c_m (H) \le \mathsf c_k (H) \le
\mathsf c (H) = m \,,
\]
and the assertion follows.

\smallskip
3. Suppose $k \in \N_{\ge 2}$ and $\mathsf{c}_{k}(H) > \mathsf{c}_{k-1}(H)$.
Let $a \in H$ with  $\min \mathsf{L}(a) \le k$ such that $\mathsf{c}(a) = \mathsf{c}_{k}(H)$.
We note that actually $\min \mathsf{L}(a) = k$, as otherwise $\mathsf{c}_{k-1}(H) \ge \mathsf{c}(a)$, a contradiction.
Let $z,z' \in \mathsf{Z}(a)$ such that $\mathsf{d}(z,z')= \mathsf{c}(a) = \mathsf{c}_{k}(H)$ and such that $z$ and $z'$ cannot be concatenated by an $N$-chain for $N < \mathsf{c}(a)$.
Let $x= \gcd(z,z')$. We note that $\min \{|x^{-1}z|, |x^{-1}z'|\} \ge k$, as otherwise $x^{-1}z$ and $x^{-1}z'$
can be concatenated by a $\mathsf{c}_{k-1}(H)$-chain, implying that $z$ and $z'$ can be concatenated by such a chain.
Thus, $\mathsf{d}(z,z')\ge k$, establishing the claim.

\smallskip
4. It suffices to show that, for every $d \in \Delta (H)$, there is a
$k \in \N$ with $2 \le k < \mathsf c (H)$ and $d \le \mathsf
c_k (H) - k$. Let $d \in \Delta (H)$. Then there is an element $a
\in H$ and factorizations $z, z' \in \mathsf Z (a)$ such that $|z'|
- |z| = d$ and $\mathsf L (a) \cap [|z|, |z'|] = \{|z|, |z'| \}$.
For $N = \min \{|z'|, \mathsf c (H)\}$, there is an $N$-chain $z =
z_0, \ldots, z_l = z'$ of factorizations from $z$ to $z'$. We may
suppose that this chain cannot be refined. This means that, for any
$i \in [1, l]$, there is no $d_i$-chain concatenating $z_{i-1}$ and
$z_i$ with $d_i < \mathsf d (z_{i-1}, z_i)$. There exists some $i
\in [1, l]$ such that $|z_{i-1}| \le |z| < |z'| \le |z_i|$, say
$z_{i-1} = x v_1 \cdot \ldots \cdot v_s$ and $z_i = x w_1 \cdot
\ldots \cdot w_t$, where $x = \gcd (z_{i-1}, z_i)$, $s, t \in \N$ and
$v_1, \ldots, v_s, w_1, \ldots, w_t \in \mathcal A (H_{\red})$. We
set $b = \pi (v_1 \cdot \ldots \cdot v_s)$, $k = \min \mathsf L (b)$
and get that
\[
2 \le k \le s < t = \max \{s, t\} = \mathsf d (z_{i-1}, z_i) =
\mathsf d (v_1 \cdot \ldots \cdot v_s, w_1 \cdot \ldots \cdot w_t)
\le N \le \mathsf c (H)  \,.
\]
Since the two factorizations $v_1 \cdot \ldots \cdot v_s$ and $w_1
\cdot \ldots \cdot w_t$ of $b$ can be concatenated by a $\mathsf c_k
(H)$-chain and since the original chain $z_0, \ldots, z_l$ cannot be
refined, it follows that $t=\mathsf d (v_1 \cdot \ldots \cdot v_s, w_1
\cdot \ldots \cdot w_t) \le \mathsf c_k (H)$. Therefore, since $|z_{i-1}| \le |z| < |z'| \le |z_i|$,  it follows
that
\[
d = |z'| - |z| \le |z_i| - |z_{i-1}| = t-s \le \mathsf c_k (H) - k
\,.
\]
Now suppose that $\mathsf c (H) \in \N$. By part 2, there is some
minimal $m \in \N$ with $\mathsf c (H) = \mathsf c_m (H)$. Since
$\mathsf c (H) > 0$, it follows that $m \ge 2$. Let $k \in \N_{\ge
2}$. If $k \ge m$, then $\mathsf c (H) = \mathsf c_m (H) = \mathsf
c_k (H)$ and $\mathsf c_k (H) -k \le \mathsf c_m (H) - m$. Thus the
assertion follows.

\smallskip
5. The inequalities are clear. Suppose that $p \in H$ is a prime
element. Let $N \in \N$ and $a \in H$ with $\mathsf c (a) \ge N$ and
$\min \mathsf L (a) \le k$. Then, for $t = k - \min \mathsf L (a)$,
we have $\mathsf L (a p^t) = t + \mathsf L (a)$, $\min \mathsf L (a
p^t) = k$ and $\mathsf c (a p^t) = \mathsf c (a) \ge N$. This
implies that
\[
\sup \{  \mathsf c (a) \mid a \in H, \
      \text {with} \ k = \min \mathsf L (a) \} \ge \sup \{ \mathsf c (a) \mid a \in H \ \text{with} \
\min \mathsf L (a) \le k \} \,,
\]
and thus equality holds in both inequalities.

\smallskip
6. Suppose that $H$ is not factorial. We start with the left
inequality. If $\mathsf L (uv) = \{2\}$ for all $u, v \in \mathcal A
(H)$, then $\daleth (H) = 0 \le \min \bigl\{\sup \Delta (H) + 2,
\mathsf c_2 (H) \bigr\}$. Let $u, v \in \mathcal A (H)$ with
$\mathsf L (uv) = \{2, d_1, \ldots , d_l\}$ with $l \in \N$ and $2 <
d_1 < \ldots < d_l$. Then $d_1 - 2 \in \Delta \bigl( \mathsf L (uv)
\bigr) \subset \Delta (H)$, and thus we get $\daleth (H) - 2 \le
\sup \Delta (H)$. Let $z' = w_1 \cdot \ldots \cdot w_{d_1} \in
\mathsf Z (uv)$ be a factorization of length $d_1$. Then, from the definition of $d_1$, we see $z = uv$ and $z'$ cannot be
concatenated by a $d$-chain with $d < d_1$. Thus $d_1 \le \mathsf c
(uv) \le \mathsf c_2 (H)$, and hence $\daleth (H) \le \mathsf c_2
(H)$.

To verify the right inequality, note that $\mathsf c_2 (H) \le
\mathsf c (H)$ follows from the definition. If $b \in H$ with
$|\mathsf Z (b)|
> 1$, then \eqref{3.0} shows that $2 + \sup \Delta \bigl(
\mathsf L (b) \bigr) \le \mathsf c (b) \le \mathsf c (H)$, and
therefore $2 + \sup \Delta (H) \le \mathsf c (H)$.
\end{proof}

\medskip
Corollary \ref{4.3} will show that, for the Krull monoids under
consideration, equality holds throughout \eqref{crucial-inequ}. Obviously, such a result is far from being true
in general. This becomes  clear from the characterization of the
catenary degree in terms of minimal relations,  recently given by S.
Chapman et al. in \cite{C-G-L-P-R06}. But we will demonstrate
this by very explicit examples which also deal with the
refinements $\mathsf c_k (H)$.

\medskip
\begin{examples} \label{3.3}~

\smallskip
{\bf 1. Numerical monoids.} The arithmetic of numerical monoids has
been studied in detail in recent years (see \cite{B-C-K-R06,
C-G-L09, A-C-H-P07a,C-G-L-M09,C-K-D-H10, Sc-Ho-Ka09a, Om10a} and the
monograph \cite{Ro-GS09}). The phenomena we are looking at here can
already be observed in the most simple case where the numerical
monoid has two generators.

Let $H = [\{d_1, d_2\}] \subset (\N_0, +)$ be a numerical monoid
generated by integers $d_1$ and $d_2$, where $1 < d_1 < d_2$ and $\gcd (d_1, d_2) =
1$. Then $\mathcal A (H)= \{d_1, d_2\}$, and $d_1 d_2$ is the
smallest element $a \in H$---with respect to the usual order in
$(\N_0, \le)$---with $|\mathsf Z (a)|
> 1$. Thus $\mathsf c_k (H) = 0$ for all $k < d_1$ (hence $\daleth (H) = 0$ if $d_1 > 2$), $\Delta (H) = \{d_2-d_1\}$ and $\mathsf
c_{d_1} (H) = d_2 = \mathsf c (H)$ (details of all this are worked
out in \cite[Example 3.1.6]{Ge-HK06a}). Thus, when $d_1 > 2$, the second two inequalities in Lemma \ref{3.2}.6 are strict.

\smallskip
{\bf 2. Finitely primary monoids.} A monoid $H$ is  called \ {\it
finitely primary}
 \ if there exist $s,\, \alpha \in \N$ with the following properties:

\begin{enumerate}
\item[]
$H$ is a submonoid of a factorial monoid $F= F^\times \time
[p_1,\ldots,p_s]$ with $s$ pairwise non-associated prime elements
$p_1, \ldots, p_s$ satisfying
\[\qquad
H \setminus H^\times \subset p_1 \cdot \ldots \cdot p_sF \quad
\text{and} \quad (p_1 \cdot \ldots \cdot p_s)^\alpha F \subset H \,.
\]
\end{enumerate}
The multiplicative monoid of every one-dimensional local noetherian
domain $R$ whose integral closure $\overline R$ is a finitely
generated $R$-module is finitely primary (\cite[Proposition
2.10.7]{Ge-HK06a}). Moreover, the monoid of invertible ideals of an
order in a Dedekind domain is a product of a free monoid and a
finite product of finitely primary monoids (see \cite[Theorem
3.7.1]{Ge-HK06a}).

Let $H$ be as above with $s \ge 2$. Then $3 \le \mathsf c (H) \le
2\alpha + 1$, $\min \mathsf L (a) \le 2 \alpha$ for all $a \in H$,
and hence $\sup \{ \mathsf c (a) \mid a \in H \
      \text {with} \ k = \min \mathsf L (a) \} = 0$ for all $k > 2
      \alpha$ (see \cite[Theorem 3.1.5]{Ge-HK06a}). This shows that
the assumption in Lemma \ref{3.2}.5 requiring the existence of a
prime element cannot be omitted. Concerning the inequalities in
Lemma \ref{3.2}.6, equality throughout can hold (as in
\cite[Examples 3.1.8]{Ge-HK06a}) but does not hold necessarily, as
the following example shows. Let $H \subset (\N_0^s, +)$, with $s\geq 3$, be the
submonoid generated by
\[
A = \{ (m,1, \ldots, 1), (1,m,1, \ldots, 1), \ldots, (1, \ldots, 1,
m) \mid m \in \N \} \,.
\]
Then $H$ is finitely primary with $A = \mathcal A (H)$ and $\daleth
(H) = 0 < \mathsf c (H)$.

\smallskip
{\bf 3. Finitely generated Krull monoids.} Let $G$ be an abelian
group and  $r, \,n \in \N_{\ge 3}$ with $n \ne r+1$.
Let \ $e_1,
\dots , e_r \in G$ be independent elements with \ $\ord (e_i) = n$
for all $i \in [1,r]$, \ $e_0 = -(e_1 + \ldots + e_r)$ and $G_0 =
\{e_0, \dots , e_r\}$. Then $\mathcal B (G_0)$ is a finitely
generated Krull monoid, $\Delta (G_0) = \{|n-r-1|\}$, $\mathsf c
(G_0) = \max \{n, r+1\}$ and
\[
0 = \daleth (H) = \mathsf c_2 (H) < 2 + \max \Delta (H) < \mathsf c
(H) \,.
\]
(see \cite[Proposition 4.1.2]{Ge-HK06a}).

\smallskip
{\bf 4. $k$-factorial monoids.} An atomic monoid $H$ is called
$k$-factorial, where $k \in \N$, if every element $a \in H$ with
$\min \mathsf L (a) \le k$ has unique factorization; $k$-factorial
and, more generally, quasi-$k$-factorial monoids and domains have been
studied in \cite{A-C-HK-Z98}. Clearly, if $H$ is $k$-factorial but
not $k+1$-factorial, then $0 = \mathsf c_k (H) < \mathsf c_{k+1}
(H)$.

\smallskip
{\bf 5. Half-factorial monoids.} An atomic monoid $H$ is called half-factorial if
$\Delta(H)= \emptyset$ (cf.~\cite[Section 1.2]{Ge-HK06a}). Then, $\daleth(H)=0$ and it follows that $\mathsf{c}_k(H)\le k$ for each $k \in \mathbb{N}$.
Thus, by Lemma \ref{3.2}.3, we get that  if $\mathsf{c}_k(H)>\mathsf{c}_{k-1}(H)$, then $\mathsf{c}_k(H) = k$.
Without additional restriction on $H$, the set $K \subset \mathbb{N}_{\ge 2}$ of all $k$ with $\mathsf{c}_k(H)>\mathsf{c}_{k-1}(H)$ can be essentially arbitrary; an obvious restriction is that it is finite for $\mathsf{c}(H)$ finite.
\end{examples}

\medskip
The arithmetic of Krull monoids is studied via transfer
homomorphisms. We recall the required terminology and collect
the results needed for the sequel.

\medskip
\begin{definition} \label{3.4}
A monoid homomorphism \ $\theta \colon H \to B$ is called a \ {\it
transfer homomorphism} \ if it has the following properties:

\smallskip

\begin{enumerate}
\item[]
\begin{enumerate}
\item[{\bf (T\,1)\,}] $B = \theta(H) B^\times$ \ and \ $\theta
^{-1} (B^\times) = H^\times$.

\smallskip

\item[{\bf (T\,2)\,}] If $u \in H$, \ $b,\,c \in B$ \ and \ $\theta
(u) = bc$, then there exist \ $v,\,w \in H$ \ such that \ $u = vw$,
\ $\theta (v) \simeq b$ \ and \ $\theta (w) \simeq c$.
\end{enumerate}\end{enumerate}
\end{definition}

\medskip
\noindent Note that the second part of (T1) means precisely that units map to units and non-units map to non-units, while the first part means $\theta$ is surjective up to units. Every transfer homomorphism $\theta$ gives rise to a
unique extension $\overline \theta \colon \mathsf Z(H) \to \mathsf
Z(B)$ satisfying
\[\qquad \quad
\overline \theta (uH^\times) = \theta (u)B^\times \quad \text{for
each} \quad u \in \mathcal A(H)\,.
\]
For $a \in H$, we denote by \ $\mathsf c (a, \theta)$ \ the smallest
$N \in \N_0 \cup \{\infty\}$ with the following property:

\smallskip

\begin{enumerate}
\item[]
If $z,\, z' \in \mathsf Z_H (a)$ and $\overline \theta (z) =
\overline \theta (z')$, then there exist some $k \in \N_0$ and
factorizations $z=z_0, \ldots, z_k=z' \in \mathsf Z_H (a)$ such that
\ $\overline \theta (z_i) = \overline \theta (z)$ and \ $\mathsf d
(z_{i-1}, z_i) \le N$ for all $i \in [1,k]$ \ (that is, $z$ and $z'$
can be concatenated by an $N$-chain in the fiber \ $\mathsf Z_H (a)
\cap \overline \theta ^{-1} (\overline \theta (z)$)\,).
\end{enumerate}

\smallskip\noindent
Then
\[
\mathsf c (H, \theta) = \sup \{\mathsf c (a, \theta) \mid a \in H \}
\in \N_0 \cup \{\infty\}
\]
denotes the {\it catenary degree in the fibres}.

\medskip
\begin{lemma} \label{3.5}
Let $\theta \colon H \to B$ be a transfer homomorphism of atomic
monoids and $\overline \theta \colon \mathsf Z (H) \to \mathsf Z(B)$
its extension to the factorization monoids.
\begin{enumerate}
\item For every \ $a \in H$, we have \ $\mathsf L_H (a) = \mathsf L_B
      \bigl( \theta (a) \bigr)$. In particular, we have  $\Delta (H) = \Delta (B)$ and $\daleth (H) = \daleth
      (B)$.

\smallskip
\item For every \ $a \in H$, we have \ $\mathsf c \bigl(\theta(a)\bigr) \le
      \mathsf c(a) \le \max \{\mathsf c\bigl( \theta(a) \bigr),\, \mathsf
      c(a,\theta)\}$.

\smallskip
\item For every $k \in \N$, we have
      \[
      \mathsf c_k (B) \le \mathsf c_k (H) \le \max \{\mathsf c_k(B), \mathsf c(H,
      \theta)\},
      \]
      and hence
      \[
      \mathsf c(B) \le \mathsf c(H) \le \max \{\mathsf c(B), \mathsf c(H, \theta)\} \,.
      \]
\end{enumerate}
\end{lemma}

\begin{proof}
1. and 2. See \cite[Theorem 3.2.5]{Ge-HK06a}.

\smallskip
3. Since, for every $a \in H$, we have $\mathsf
L (a) = \mathsf L \bigl( \theta (a) \bigr)$, it follows that $\min \mathsf L
(a) = \min \mathsf L \bigl( \theta (a) \bigr)$, and thus parts 1 and 2 imply
both inequalities.
\end{proof}

\medskip
\begin{lemma} \label{3.6}
Let $H$ be a Krull monoid, $\varphi \colon H \to F = \mathcal F (P)$
a cofinal divisor homomorphism, $G = \mathcal C (\varphi)$ its class
group, and $G_P \subset G$ the set of classes containing prime
divisors. Let $\widetilde{\boldsymbol \beta} \colon F \to \mathcal F
(G_P)$ denoted the unique homomorphism defined by
$\widetilde{\boldsymbol \beta} (p) = [p]$ for all $p \in P$.
\begin{enumerate}
\item The homomorphism $\boldsymbol \beta = \widetilde{\boldsymbol \beta} \circ \varphi \colon H \to \mathcal B
      (G_P)$ is a transfer homomorphism with $\mathsf c (H,
      \boldsymbol \beta) \le 2$.

\smallskip
\item For every $k \in \N$, we have
      \[
      \mathsf c_k (G_P) \le \mathsf c_k (H) \le \max \{\mathsf c_k(G_P), 2\},
      \]
      and hence
      \[
      \mathsf c(G_P) \le \mathsf c(H) \le \max \{\mathsf c(G_P), 2\} \,.
      \]

\smallskip
\item $\daleth (H) = \daleth (G_P) \le \mathsf D (G_P)$.
\end{enumerate}
\end{lemma}

\begin{proof}
1. This follows from  \cite[Theorem 3.4.10]{Ge-HK06a}.

\smallskip
2. This follows from part 1 and Lemma \ref{3.5}.

\smallskip
3. Since $\boldsymbol \beta$ is a transfer homomorphism, we have
$\daleth (H) = \daleth (G_P)$ by Lemma \ref{3.5}. In order to show that $\daleth (G_P)
\le \mathsf D (G_P)$, let $U_1, U_2 \in \mathcal A (G_P)$. If
$\mathsf D (G_P) = 1$, then $G_P = \{0\}$, $U=V=0$ and $\daleth
(G_P) = 0$. Suppose that $\mathsf D (G_P) \ge 2$  and consider a
factorization
 $U_1 U_2 = W_1 \cdot \ldots \cdot W_s$, where $s \in \N$ and $W_1, \ldots, W_s \in \mathcal A (G_P)$. It suffices to show that
$s \le \mathsf D (G_P)$. For $i \in [1, s]$, we set $W_i =
W_i^{(1)}W_i^{(2)}$ with $W_i^{(1)}, W_i^{(2)} \in \mathcal F (G_P)$
such that $U_1 = W_1^{(1)} \cdot \ldots \cdot W_s^{(1)}$ and $U_2 =
W_1^{(2)} \cdot \ldots \cdot W_s^{(2)}$. If  there are $i \in [1,s]$
and $j \in [1,2]$, say $i=j=1$, such that $W_i^{(j)}=W_1^{(1)} = 1$, then $W_1
= W_1^{(2)} \t U_2$; hence $W_1 = U_2$, $W_2 = U_1$ and $s = 2 \le
\mathsf D (G_P)$. Otherwise, we have $W_1^{(j)}, \ldots, W_s^{(j)}
\in \mathcal F (G_P) \setminus \{1\}$, and hence $s \le
\sum_{i=1}^s|W_i^{(j)}| = |U_j| \le \mathsf D (G_P)$.
\end{proof}

\bigskip
\section{A structural result for the catenary degree} \label{4}
\bigskip

In Theorem \ref{4.2} we obtain a structural result for the catenary degree.
Since it is relevant for the discussion of this result, we start with a technical result.

\medskip
\begin{proposition} \label{4.1}
Let $G$ be an abelian group.
\begin{enumerate}
\item  Let $G_0 = \{e_0, \ldots, e_r, -e_0, \ldots , -e_r\} \subset G$
       be a subset with $e_1, \ldots, e_r \in G$ independent
       and $e_0 = k_1 e_1 + \ldots + k_r e_r$, where  $k_i \in \N$ and $2k_i
       \le \ord (e_i)$ for all $i \in [1, r]$. If $\sum_{i=1}^r k_i\neq 1$, then $\daleth (G_0)
       \ge k_1 + \ldots + k_r+1$.

\smallskip

\item Let $G_0=\{-e,e\} \subset G$ be a subset with $3 \le \ord (e)< \infty$. Then
$\daleth(G_0)\ge \ord(e)$.

\smallskip
\item Let $G = C_{n_1} \oplus \ldots \oplus C_{n_r}$ with $|G| \ge 3$ and $1 < n_1 \t
      \ldots \t n_r$, and let $(e_1, \ldots, e_r)$ be a basis of $G$ with $\ord (e_i)
           = n_i$ for all $i \in [1, r]$. If  $\{e_0, \ldots,
           e_r, -e_0, \ldots , -e_r\} \subset G_0 \subset G$, where $e_0 = \sum_{i=1}^r
           \lfloor \frac{n_i}{2} \rfloor e_i$, then $\daleth (G_0) \ge
      \max\{n_r,\,1+\sum_{i=1}^{r}\lfloor\frac{n_i}{2}\rfloor\}$.
\end{enumerate}
\end{proposition}

\begin{proof}
1.
If
\[
A = e_0(-e_0) \prod_{i=1}^r e_i^{k_i} (-e_i)^{k_i} \\,
\]
then $\mathsf L (A) = \{2, k_1+ \ldots + k_r+1\}$ (see \cite[Lemma
6.4.1]{Ge-HK06a}). Thus, if $\sum_{i=1}^r k_i \neq 1$, the assertion follows by definition of
$\daleth (G_0)$.

\smallskip
2. Let $n = \ord(e)$. Since $\mathsf L \bigl( (-e)^{n} e^{n} \bigr) = \{2, n\}$, we get $\daleth
(G_0) \ge n$.

\smallskip
3. Clear, by parts 1 and 2.
\end{proof}

\medskip
\begin{theorem} \label{4.2}
Let $H$ be a Krull monoid, $\varphi \colon H \to F = \mathcal F (P)$
a cofinal divisor homomorphism, $G = \mathcal C (\varphi)$ its class
group, and $G_P \subset G$ the set of classes containing prime
divisors. Then \be \label{whoopcou}\mathsf c(H) \leq \max \Big\{
\Big\lfloor\frac{1}{2} \mathsf D(G_P)+1 \Big\rfloor,\, \daleth (G_P)
\Big\}. \ee
\end{theorem}

\begin{proof}
By Lemma \ref{3.6}, we have $\mathsf c (H) \le \max \{\mathsf c
(G_P), 2\}$. If $\mathsf D (G_P) = 1$, then $G_P = \{0\}$, $G =
[G_P] = \{0\}$, $H = F$ and $\mathsf c (H) = 0$. Thus we may suppose
that $2 \le \mathsf D (G_P) < \infty$, and it is sufficient to show
that
\[
\mathsf c (G_P) \le d_0, \quad \text{where} \quad d_0 = \max \Big\{
\Bigl\lfloor\frac{1}{2} \mathsf D(G_P)+1 \Big\rfloor,\, \daleth
(G_P) \Big\} \,.
\]
So we have to verify that, for $A\in \mathcal B(G_P^{\bullet})$ and
$z,\,z'\in \mathsf Z(A)$,  there is a $d_0$-chain of factorizations
between $z$ and $z'$. Assuming this is false, consider a counter
example $A\in \mathcal B(G_P^{\bullet})$ such that $|A|$ is minimal,
and for this $A$, consider a pair of factorizations $z,\,z'\in
\mathsf Z(A)$ for which no $d_0$-chain between $z$ and $z'$ exists
such that $|z|+|z'|$ is maximal (note $|A|$ is a trivial upper bound for the length
of a factorization of $A$).

Note we may assume \be\label{rewind}\max\{|z|,\,|z'|\}\geq d_0+1\geq
\frac{1}{2} \mathsf D (G_P)+\frac32,\ee else the chain $z,z'$ is a
$d_0$-chain between $z$ and $z'$, as desired. We continue with the
following assertion.

\begin{enumerate}
\item[{\bf A.}\,]
Let
\[
y = U_1 \cdot \ldots \cdot U_r\in \mathsf Z(A) \quad \text{and}
\quad y' = V_1 \cdot \ldots \cdot V_s\in \mathsf Z(A) \,, \quad
\text{where} \quad U_i,\,V_j\in \mathcal A(G_P)\,,
\]
be two factorizations of $A$ with $V_{j_1}|U_1 \cdot \ldots \cdot
U_rU_{j_2}^{-1}$, for some $j_1\in [1,s]$ and $j_2\in [1,r]$. Then
there is a  $d_0$-chain of factorizations of $A$ between $y$ and
$y'$.
\end{enumerate}

\smallskip
{\it Proof of \,{\bf A}}.\, We may
assume $j_1=1$, $j_2=r$, and we obtain a factorization
\[
U_1 \cdot \ldots \cdot  U_{r - 1} = V_1 W_1  \cdot \ldots \cdot W_t
\,,
\]
where $W_1, \ldots , W_t \in \mathcal A (G_P)$. By the minimality of
$|A|$, there is a $d_0$-chain of factorizations $y_0, \ldots, y_k$
between $y_0 = U_1 \cdot \ldots \cdot U_{r - 1}$ and $y_k = V_1 W_1
\cdot \ldots \cdot W_t$, and there is a $d_0$-chain of
factorizations $z_0, \ldots, z_l$ between $z_0 = W_1 \cdot \ldots
\cdot W_t U_r$ and $z_l = V_2 \cdot \ldots \cdot V_s$. Then
\[
y= y_0  U_r, \, y_1 U_r, \ldots, y_k U_r =  V_1 z_0, \, V_1 z_1,
\ldots, V_1  z_l  = y'
\]
is a $d_0$-chain between $y$ and $y'$. \qed

\smallskip

We set $z = U_1 \cdot \ldots \cdot U_r$ and $z' = V_1 \cdot \ldots
\cdot V_s$, where all $U_i,\,V_j \in \mathcal A (G_P)$, and without
loss of generality we assume that $r \ge s$. Then, in view of
\eqref{rewind} and $\mathsf D(G_P)\geq 2$, it follows that \be\label{hant}r \ge d_0+1\geq
\frac{1}{2} \mathsf D (G_P)+\frac32 > 2 .\ee Clearly, $s=1$ would
imply $r=1$, and thus we get $s \ge  2$.

Suppose $\max \mathsf L(V_{1}V_{2})\geq 3$.
Then, by definition of $\daleth (G_P)$, there exists $y\in \mathsf Z
(V_1V_2)$ with $3 \le |y|\le \daleth (G_P)$ and
\be\label{allhallowe}
\mathsf d(z',yV_3\cdot\ldots \cdot V_s) = \mathsf d (V_1V_2, y) = |y| \le \daleth (G_P) \,.
\ee
But, since $|z|+|y V_3\cdot\ldots \cdot V_s|>|z|+|z'|$, it follows,
from the maximality of $|z|+|z'|$, that there is a $d_0$-chain of
factorizations between $y V_3\cdot\ldots \cdot V_s$ and $z$, and
thus, in view of \eqref{allhallowe}, a $d_0$-chain concatenating $z'$ and $z$, a
contradiction. So we may instead assume  $\max \mathsf L(V_{1}V_{2})=2$.

As a result, if $s=2$, then $V_1V_2 = A $ and $\mathsf L (A) = \{2\}$, contradicting  $2 < r \in \mathsf L (A)$ (cf.~\eqref{hant}).
Therefore we have $s \ge 3$.

We set $V_1 = V_1^{(1)}\cdot \ldots \cdot V_1^{(r)}$ and $V_2 =
V_2^{(1)} \cdot \ldots \cdot V_2^{(r)}$, where
$V_1^{(j)}V_2^{(j)}|U_j$ for all $j \in [1, r]$. In view of
\textbf{A}, we see that each $V_1^{(i)}$ and $V_2^{(j)}$ is
nontrivial. Thus \eqref{hant} implies \be\label{gant}|V_1V_2|\geq
2r\geq \mathsf D (G_P)+3.\ee By the pigeonhole principle and in view
of \eqref{hant}, there exists some  $j \in [1, r]$, say  $j=r$, such
that
\[
|V_1^{(r)}V_2^{(r)}|\leq \frac{1}{r}|V_1V_2|\leq \frac{2 \mathsf D
(G_P)}{r}<4 \,.
\]
As a result, it follows in view of \eqref{gant} that
\be\label{hhant}|V_1^{(1)} \cdot \ldots \cdot V_1^{(r-1)}V_2^{(1)}
\cdot \ldots \cdot V_2^{(r-1)}|\geq |V_1V_2|-3\geq \mathsf D
(G_P).\ee Thus there exists a $W_1 \in \mathcal A (G_P)$ such that
$W_1 \t V_1^{(1)} \cdot \ldots \cdot V_1^{(r-1)}V_2^{(1)} \cdot
\ldots \cdot V_2^{(r-1)}$.

Let $V_1V_2=W_1 \cdot \ldots \cdot W_t$, where $W_2, \ldots, W_t \in
\mathcal A(G_P)$. Since $s\geq 3$, we have $|V_1V_2|<|A|$. Thus, by
the minimality of $|A|$, there is a $d_0$-chain of factorizations
between $V_1  V_2$ and $W_1 \cdot \ldots \cdot W_t$, and thus one
between $z'=(V_1  V_2) V_3\cdot\ldots\cdot V_s$ and $(W_1 \cdot
\ldots \cdot W_t) V_3\cdot\ldots \cdot V_s$ as well. From the
definitions of the $V_i^{(j)}$ and $W_1$, we have $W_1 \t U_1 \cdot
\ldots \cdot U_{r-1}$. Thus by \textbf{A} there is a $d_0$-chain of
factorizations between $W_1 \cdot \ldots \cdot W_t V_3\cdot
\ldots\cdot V_s$ and $z=U_1\cdot\ldots\cdot U_r$. Concatenating
these two chains gives a $d_0$-chain of factorizations between $z'$
and $z$, completing the proof.
\end{proof}

\medskip
\begin{corollary} \label{4.3}
Let $H$ be a Krull monoid, $\varphi \colon H \to F = \mathcal F (P)$
a cofinal divisor homomorphism, $G = \mathcal C (\varphi) \cong C_{n_1}
\oplus \ldots \oplus C_{n_r}$ its class group, where $1<n_1|\ldots |n_r$ and $|G|\ge 3$, and  $G_P \subset G$  the set of all classes containing
prime divisors. Suppose that the following two conditions hold{\rm
\,:}
\begin{itemize}
\item[(a)] $\big\lfloor\frac{1}{2}\mathsf D (G_P)+1 \big\rfloor \le
           \max \big\{n_r,\,1+\sum_{i=1}^{r}\lfloor\frac{n_i}{2}\rfloor \big\}$.

\smallskip
\item[(b)] There is a basis $(e_1, \ldots , e_r)$ of $G$ with $\ord (e_i)
           = n_i$, for all $i \in [1, r]$, such that \newline $\{e_0, \ldots,
           e_r, -e_0, \ldots , -e_r\} \subset G_P$, where $e_0 =  \sum_{i=1}^r
           \lfloor \frac{n_i}{2} \rfloor e_i$.
\end{itemize}
Then
\[
\daleth (H) = 2 + \max \Delta(H) =  \mathsf c_2(H) = \mathsf c (H)
\,.
\]
\end{corollary}

Before giving the proof of the above corollary, we analyze the result
and its assumptions.

\medskip
\begin{remark} \label{4.4}
Let all notation be as in Corollary \ref{4.3}.

\smallskip
1. Note that
\[
1+\sum_{i=1}^{r} \Big\lfloor\frac{n_i}{2} \Big\rfloor =
1+\frac{\mathsf{r}_2(G)+ \mathsf d^* (G)}{2} \,,
\]
where $\mathsf{r}_2(G)$ denotes the $2$-rank of $G$, i.e., the number of even $n_i$s.
Thus, if $\mathsf D(G)=\mathsf d^*(G)+1$ (see the comments after \eqref{soitis} for some groups fulfilling this), then
\[
\Big\lfloor\frac{1}{2}\mathsf D(G)+1 \Big\rfloor
\le 1+\sum_{i=1}^{r}\Big\lfloor\frac{n_i}{2} \Big\rfloor \,,
\]
and hence Condition $(a)$ holds. Not much is known about groups $G$
with $\mathsf D (G) > \mathsf d^* (G) + 1$ (see \cite{Ge-Sc92},
\cite[Theorem 3.3]{Ga-Ge99}). Note that groups of odd order with $\mathsf D (G) > \mathsf d^* (G) + 1$ yield examples of groups for which (a) fails,
yet the simplest example of such a group we were able to find in the literature already has rank $8$ (see \cite[Theorem 5]{Ge-Sc92}).

\smallskip
2.  In Examples \ref{3.3}, we pointed out that some assumption on
$G_P$ is needed in order to obtain the result $\daleth (H) = \mathsf c
(H)$. Clearly, Condition $(b)$ holds if every class contains a prime
divisor. But since there are relevant Krull monoids with $G_P \ne G$
(for examples arising in the analytic theory of Krull monoids, we
refer to \cite{Ge-Ka92, IR99a, IR99b}), we formulated our
requirements on $G_P$ as weak as possible, and we discuss two
natural settings which enforce parts of Conditions $(b)$ even if
$G_P \ne G$.

$(i)$ A Dedekind domain $R$ is a quadratic extension of a principal
ideal domain $R'$ if $R' \subset R$ is a subring and $R$ is a free
$R'$-module of rank $2$. If $R$ is such a Dedekind domain, $G$ its
class group, and $G_P \subset G$ the set of classes containing prime
divisors, then $G_P = -G_P$ and $[G_P] = G$. By a result of
Leedham-Green \cite{Le72a}, there exists, for every abelian group
$G$, a Dedekind domain $R$ which is a quadratic extension of a
principal ideal domain and whose class group is isomorphic to $G$.

$(ii)$ If $G_P \subset G$ are as in
Corollary \ref{4.3},
then  $G_P$ is a
generating set of $G$, and if $G \cong C_{p^k}^r$, where $p \in
\mathbb P$ and  $k, r \in \N$, then $G_P$ contains a basis by
\cite[Lemma A.7]{Ge-HK06a}.

\smallskip
3.  Corollary \ref{4.3} tells us that the catenary degree $\mathsf c
(H)$ occurs as a distance of two factorizations of the following
form
\[
a = u_1 u_2 = v_1 \cdot \ldots \cdot v_{\mathsf c (H)} \,,
\]
where $u_1, u_2, v_1, \ldots, v_{\mathsf c (H)} \in \mathcal A (H)$
and $a$ has no factorization of length $j \in [3, \mathsf c (H)-1]$.
Of course,  the catenary degree may also occur as a distance between
factorizations which are not of the above form. In general, there
are even elements $a$ and integers $k \ge 3$ such that
\begin{equation}
\mathsf c (a) = \mathsf c (H) \,, \ \min \mathsf L (a) = k \quad
\text{and} \quad \mathsf c (b) < \mathsf c (a) \label{4.0}
\end{equation}
for all proper divisors $b$ of $a$. We provide a simple, explicit
example.

Let $G = C_3 \oplus C_3$, $(e_1, e_2)$ be a basis of $G$ and $e_0 =
-e_1-e_2$. For $i \in [0,2]$, let $U_i = e_i^3$ and let $V = e_0e_1e_2$.
Then $A = V^3 \in \mathcal B (G)$, $\mathsf Z (A) = \{U_0 U_1 U_2,
V^3 \}$, $\mathsf c (A) = 3 = \mathsf c (G)$ (see Corollary
\ref{5.5}) and $\mathsf c (B) = 0$ for all proper zero-sum
subsequences $B$ of $A$.

\smallskip
4. Let $\boldsymbol \beta \colon H \to \mathcal B (G_P)$ be as in
Lemma \ref{3.6}. Clearly, if $a \in H$ is such that $\mathsf c (a) =
\mathsf c (H)$, then, using the notation of Remark \ref{4.4}.3, $a, \boldsymbol \beta (a)$, $u_1$, $u_2$, $\boldsymbol \beta (u_1)$ and
$\boldsymbol \beta (u_2)$ must be highly structured. On the opposite
side of the spectrum, there is the following result: if $\supp \big(
\boldsymbol \beta (a) \big) \cup \{0\}$ is a subgroup of $G$, then
$\mathsf c (a) \le 3$ (see \cite[Theorem 7.6.8]{Ge-HK06a}),  while \eqref{3.0} shows $\mathsf c (a) \ge 3$ whenever
$|\mathsf L (a)| > 1$.

\smallskip
5. If $H$ is factorial, in particular if $|G|=1$, then $\daleth(H)=\mathsf c_2(H)=\mathsf c(H)=0$ and $2+\max \Delta(H)=2$.
If $H$ is not factorial and $|G|=2$, then $\daleth(H)=0$ and $\mathsf c_2(H) = \mathsf c(H) = 2+\max \Delta(H)=2$.
\end{remark}

\medskip
\begin{proof}[Proof of Corollary \ref{4.3}]
Lemma \ref{3.2}.6  and Theorem \ref{4.2}  imply that
\[
  \begin{split}    \daleth (H) & \le \min\{2+\max \Delta (H) ,\,\mathsf c_2(H)\} \le \max\{2+\max \Delta (H) ,\,\mathsf c_2(H)\}  \\
   & \le \mathsf c (H)  \le \max \Big\{ \Big\lfloor\frac{1}{2}\mathsf D(G_P)+1 \Big\rfloor,\, \daleth (G_P) \Big\} \,.
\end{split}
\]
By assumption and by Proposition \ref{4.1} and Lemma \ref{3.6}.3, it follows that
\[
\Big\lfloor\frac{1}{2}\mathsf D(G_P)+1 \Big\rfloor \le \max
\Big\{n_r,\,1+\sum_{i=1}^{r} \Big\lfloor\frac{n_i}{2} \Big\rfloor
\Big\} \le \daleth (G_P) = \daleth (H) \,,
\]
and thus, in the above chain of inequalities, we indeed have equality
throughout.
\end{proof}

\medskip
\begin{corollary} \label{4.6}
Let $H$ be a Krull monoid, $\varphi \colon H \to F = \mathcal F (P)$
a cofinal divisor homomorphism, $G = \mathcal C (\varphi)$ its class
group,  $G_P \subset G$ the set of classes containing prime
divisors, and suppose that $3 \le \mathsf D (G_P) < \infty$.
\begin{enumerate}
\item We have \ $\mathsf c (H) = \mathsf D (G_P)$ \ if and only if \ $\daleth (H) = \mathsf D (G_P)$.

\smallskip
\item If \ $\mathsf c (H) = \mathsf D (G)$, then \ $\mathsf D (G_P) = \mathsf D (G)$ \ and \ $G$ is either cyclic or an elementary $2$-group.
If $G_P=-G_P$, then the converse implication holds as well.
\end{enumerate}
\end{corollary}

\begin{proof}
1. By Theorem \ref{4.2}, \eqref{crucial-inequ} and Lemma \ref{3.6}.3, we have
\be\label{wraiter}
\daleth (H) = \daleth (G_P) \le \mathsf c (H) \le  \max \Big\{
\Big\lfloor\frac{1}{2}\mathsf D(G_P)+1 \Big\rfloor,\, \daleth (G_P)
\Big\} \le \mathsf D (G_P) \,,
\ee
which we will also use for part 2. In view of $3\leq \mathsf D(G_P)<\infty$, we have $\lfloor \frac12 \mathsf D(G_P)+1\rfloor <\mathsf D(G_P)$. Thus the assertion now directly follows from \eqref{wraiter}.

\smallskip
2. We use that $[G_P] = G$. Furthermore, if $\mathsf D (G_P) =
\mathsf D (G)$, it follows that $\Sigma (S) = G^{\bullet}$
      for all zero-sum free sequences $S \in \mathcal F (G_P)$ with $|S| =
      \mathsf D (G_P)-1$ (see \cite[Proposition 5.1.4]{Ge-HK06a}).
Obviously, this implies that $\langle \supp (U) \rangle = G$ for all
$U \in \mathcal A (G_P)$ with $|U| = \mathsf D (G_P)$.

Suppose that $\mathsf c (H) = \mathsf D (G)$. Since
$\mathsf c (H) \le \mathsf D (G_P) \le \mathsf D (G)$ (in view of \eqref{wraiter}), it follows
that $\mathsf D (G_P) = \mathsf D (G)$, and part 1 implies that
$\daleth (H) = \mathsf D (G_P)$. Thus there exist $U, V \in \mathcal
A (G_P)$ such that $\{2, \mathsf D (G) \} \subset \mathsf L  (U V
)$, and \cite[Proposition 6.6.1]{Ge-HK06a} implies that $V = -U$ and
$\mathsf L \bigl( (-U) U \bigr) = \{2, \mathsf D (G) \}$ (since $\max \mathsf L((-U)U)\leq \frac{|(-U)U|}{2}\leq \mathsf D(G)$).

Assume to the contrary that $G$ is neither cyclic nor an elementary
$2$-group. We show that there exists some $W \in \mathcal A (G_P)$
such that $W \t (-U) U$ and $2 < |W| < \mathsf D (G)$. Clearly, $W$
gives rise to a factorization $(-U)U = W W_2 \cdot \ldots \cdot W_k$
with $W_2, \ldots, W_k \in \mathcal A (G_P)$ and $2 < k < \mathsf D
(G)$, a contradiction to $\mathsf L \bigl( (-U) U \bigr) = \{2, \mathsf D (G) \}$.

Since $\langle \supp (U) \rangle = G$ (as noted above) is not an elementary
$2$-group, there exists some $g_0 \in \supp(U)$ with $\ord (g_0) >
2$, say $U = g_0^m g_1 \cdot \ldots \cdot g_l$ with $g_0 \not\in \{
g_1, \dots, g_l \}$. Since $G = \langle \supp (U) \rangle$ is not
cyclic, it follows that $l \ge 2$. Let \ $W' = (- g_0)^m g_1 \cdot
\ldots \cdot g_l$. Then $W' \t U(-U)$ and $|W'| = \mathsf  D (G)$.
Hence there exists some $W \in \mathcal A(G_P)$ with $W \t W'$, and
we proceed to show that $2<|W | < \mathsf D(G)$, which will complete the proof. Since $U \in \mathcal A
(G_P)$, we have $W \nmid g_1 \cdot \ldots \cdot g_l$, and thus $-g_0
\t W$. Since $g_0 \notin \{ g_1, \dots, g_l \}$ and $g_0 \ne -g_0$,
it follows that $W \ne g_0(-g_0)$, and thus $|W| > 2$.

Assume to the contrary that $|W| = \mathsf  D (G)$. Then $W = W'$,
and $\sigma(U) = \sigma (W')=0$ implies $2 m g_0 = 0$, and thus $m >
1$. We consider the sequence $S = g_0^m g_1 \cdot \ldots \cdot g_{l
- 1}$. Since $1<m<\ord (g_0)$ and $2mg_0=0$, it follows that $$0\neq (m+1)g_0.$$ Since $S$ is zero-sum free of length $|S| = \mathsf D
(G)-1$, we have $\Sigma (S) = G^{\bullet}$, and thus $0\neq (m + 1) g_0 \in
\Sigma (S)$, say
\[
(m + 1) g_0 = s g_0 + \sum_{i \in I} g_i \quad \text{with} \quad s
\in [0,m] \quad \text{and} \quad I \subset [1,  l - 1]\,.
\]
If  $s = 0$, then
\[
0 = 2 m g_0 = (m - 1) g_0 + \sum_{i \in I} g_i \ \in \Sigma (S) ,
\]
a contradiction. If $s \ge 1$, then it follows that
\[
T =  (- g_0)^{m + 1 - s} \prod_{i \in I} g_i
\]
is a proper zero-sum subsequence of  $W$, a contradiction to  $W \in
\mathcal A (G_P)$.

Suppose that $G_P=-G_P$ and $\mathsf D (G_P) = \mathsf D (G)$. Recall the comments after \eqref{soitis} concerning the value of $\mathsf D(G)$. First, we let $G$ be an elementary $2$-group. Then there is a $U = e_0e_1 \cdot
\ldots \cdot e_r \in \mathcal A (G_P)$ with $|U| = \mathsf D (G) =
r+1$. Thus, since $\langle \supp(U)\rangle =G$, and since a basis of an elementary $2$-group is just a minimal (by inclusion) generating set, it follows that  $G_P$ contains the basis (say) $(e_1, \ldots, e_r)$ of $G$,
and Proposition \ref{4.1} and Lemma \ref{3.6}.3 imply that $\daleth(H)=\daleth (G_P) = \mathsf D
(G_P) = \mathsf D (G)=r+1$, whence $\mathsf c(H)=\mathsf D(G)$ follows from part 1. Second,  let $G$ be cyclic. If $U \in
\mathcal A (G_P)$ with $|U| = \mathsf D (G_P) = \mathsf D (G)$, then
$|U| = |G|$ and \cite[Theorem 5.1.10]{Ge-HK06a} implies that $U =
g^{|G|}$ for some $g \in G_P$ with $\ord (g) = |G|$. Hence $\mathsf
L \bigl( (-U)U \bigr) = \{2, |G|\}$, and now it follows from Lemma \ref{3.6}.3 that $|G| = \mathsf
D (G_P) = \daleth (G_P) = \daleth (H)$, whence part 1 once more shows $\mathsf c(H)=\mathsf D(G)=\mathsf D(G_P)$.
\end{proof}

\bigskip
\section{An upper bound for the catenary degree} \label{5}
\bigskip

We apply our structural result on the catenary degree (Theorem \ref{4.2})
to obtain a new upper bound on the catenary degree (see Theorem \ref{5.4}) and
a characterization result for Krull monoids with small catenary degree (see Corollary \ref{5.6}).
We start with some technical results.

\medskip
\begin{lemma} \label{5.1}
Let $G$ be an abelian group and  let $U,\,V\in \mathcal F
(G^{\bullet})$. Suppose that either $U, V \in \mathcal A (G)$ or
that $U$ and $V$ are zero-sum free with $\sigma (UV) = 0$. Then
$\max \mathsf L(UV) \le \min \{|U|,\,|V|\}$. Moreover, if $\max
\mathsf L(UV)=|U|\geq 3$, then $-\supp(U)\subset \Sigma(V)$.
\end{lemma}

\begin{proof}
Let $UV = W_1 \cdot \ldots \cdot W_m$, where $m = \max \mathsf L
(UV)$ and $W_1, \ldots, W_m \in \mathcal A (G)$.
 Let  $U = U_1 \cdot \ldots \cdot U_m$
and $V = V_1 \cdot \ldots \cdot V_m$ with $W_i = U_i V_i$ for $i \in
[1, m]$. If $U_i \ne 1$ and $V_i \ne 1$ for all $i \in [1, m]$, then
$m \le |U_1| + \ldots + |U_m| = |U|$ and likewise $m \le |V|$.
Moreover, if equality holds in the first bound, then $|U_i|=1$ for
$i\in [1,m]$, in which case each $V_i|V$ is a subsequence of $V$
with $\sigma(V_i)=-\sigma(U_i)\in -\supp(U)$; since
$\bigcup_{i=1}^{m}\{ \sigma(U_i) \} = \supp(U)$, this means
$-\supp(U)\subset \Sigma(V)$.

On the other hand, if there is some $j \in [1, m]$ such that $U_j = 1$
or $V_j = 1$, say $U_1 = 1$, then, since $V$ contains no proper, nontrivial zero-sum subsequence, it follows that  $W_1 = V_1 = V$, which, since $U$ contains no proper, nontrivial zero-sum subsequence, implies $W_2
= U$. Hence, since $U,\,V\in \mathcal F(G^{\bullet})$ with $\sigma(U)=\sigma(W_2)=0=\sigma(W_1)=\sigma(V)$ implies $|U|,\,|V|\geq 2$, we see that $m = 2 \le \min \{|U|,\,|V|\}$.
\end{proof}

\medskip
\begin{lemma}\label{5.2}
Let $G$ be an abelian group,  $K \subset  G$  a finite cyclic
subgroup, and let $U,\,V\in \mathcal A(G)$ with $\max \mathsf L (UV)
\geq 3$. If $\sum_{g\in K} \mathsf v_g(UV)\geq |K|+1$ and there
exists a nonzero $g_0\in K$ such that $\mathsf v_{g_0}(U)>0$ and
$\mathsf v_{-g_0}(V)>0$, then $\mathsf L (UV)\cap [3,|K|] \ne
\emptyset$.
\end{lemma}

\begin{proof}
Note $U,\,V\in \mathcal A(G)$ and $\max \mathsf L (UV) \geq 3$ imply
$0\notin \supp(UV)$.
Moreover, note that if $\supp(U)\subset K$, then  Lemma \ref{5.1} implies that $\max
\mathsf L (UV)\le |U|\le \mathsf D(K)= |K|$ (recall the comments after \eqref{soitis}), whence the assumption
$\max \mathsf L (UV) \geq 3$ completes the proof. Therefore we may
assume $\supp(U)\not\subset K$, and likewise that
$\supp(V)\not\subset K$.

We factor $U=U_0U'$ and $V=V_0V'$ where $U_0$
and $V_0$ are subsequences of terms from $K$ such that there exists some non-zero $g_0 \in K$ with $g_0\t U_0 $ and $(-g_0)\t V_0$, and $|U_0|+|V_0|=|K|+1$.
Note that by the assumption made above, both $U_0$ and $V_0$ are
proper subsequences of $U$ and $V$, respectively, and thus they are zero-sum free.

Let $U_0=g_0U_0'$ and $V_0 = (-g_0)V_0'$.
Since $U_0'$ and $V_0'$ are both zero-sum free,
we get (cf., e.g., \cite[Proposition 5.1.4.4]{Ge-HK06a}) that $|\{0\} \cup \Sigma(U_0')| \ge |U_0'| + 1 = |U_0|$
and $|\{0\} \cup \Sigma(V_0')| \ge |V_0'| + 1 = |V_0|$.
Since these sets are both subsets of $K$, the pigeonhole principle implies that
\begin{equation}
\label{eq_nonempty}
\bigl(g_0 + \bigl( \{0\} \cup \Sigma(U_0') \bigr) \bigr) \cap \bigl( \{0\} \cup \Sigma(V_0') \bigr) \neq \emptyset .
\end{equation}
Let $U_0''$ and $V_0''$ denote (possibly trivial) subsequences of $U_0'$ and $V_0'$, respectively, such that
$\sigma(V_0'')= g_0+\sigma(U_0'')=\sigma(g_0U_0'')$, whose existence is guaranteed by \eqref{eq_nonempty}.

We set $W_1= (g_0 U_0'')^{-1} U V_0''$ and $W_2= V_0''^{-1} V (g_0U_0'')$.
Then, $UV=W_1W_2$, and $W_1$ and $W_2$ are nontrivial zero-sum sequence; more precisely, $(-g_0)g_0 \t W_2$ is a proper
zero-sum subsequence (recall that by assumption $U_0$ and $V_0$ are proper subsequences of $U$ and $V$, respectively).
Since $\mathsf{L}(W_1)+ \mathsf{L}(W_2) \subset \mathsf{L}(UV)$, and since by the above assertion
$\min \mathsf{L}(W_1) \ge 1 $ and $ \min \mathsf{L}(W_2) \ge 2$, it suffices to assert that
$\max \mathsf{L}(W_1) + \max \mathsf{L}(W_2) \le |K|$.
Since, by Lemma \ref{5.1}, we have $\max \mathsf{L}(W_1)\le |V_0''|\le |V_0|-1$ and
$\max \mathsf{L}(W_2) \le |g_0U_0''|\le |U_0|$, and  since by assumption $|U_0| + |V_0|= |K|+1$, this is the case.
\end{proof}

\medskip
\begin{lemma} \label{5.3}
Let $t \in \N$ and $\alpha, \alpha_1, \ldots, \alpha_t \in \R$ with
$\alpha_1 \ge \ldots \ge \alpha_t \ge 0$ and $\Sum{i=1}{t}\alpha_i\geq \alpha\geq 0$. Then
\[
\prod_{i=1}^t (1 + x_i)  \qquad \text{is minimal}
\]
over all $(x_1,\ldots,x_t) \in \R^t$ with $0 \le x_i \le \alpha_i$
and $\sum_{i=1}^t x_i = \alpha$ if
\[
x_i = \alpha_i \text{ for each } i \in [1,s] \quad \text{and} \quad x_i = 0 \quad \text{ for each } i \in [s+2,t]
\]
where $s \in [0, t]$ is maximal with $\sum_{i=1}^s \alpha_i \le \alpha$.
\end{lemma}

\begin{proof}
This is a simple calculus problem; for completeness, we include a short proof. We may assume $\alpha\neq 0$. By compactness and continuity, the existence of a minimum is clear.
Let $\overline{x}=(x_1, \dots, x_t)$ be a point where the minimum is attained.
We note that for $x,y \in \mathbb{R}$ with $x\ge y \ge 0$ we have
\begin{equation}\label{eq_decrease}
(1+x+\varepsilon)(1+y-\varepsilon)< (1+x)(1+y)
\end{equation}
for each  $\varepsilon > 0$.
Thus, it follows that $x_{i} \notin \{0,\alpha_{i}\}$ for at most one $i\in [1,t]$; if such an $i$ exists we denote it by $i_0$, otherwise we denote by $i_0$ the maximal $i \in [1,t]$ with $x_i \neq 0$.
Suppose that for $\overline{x}$ the value of $\alpha_{i_0}$ is maximal among all points where the minimum is attained.
We observe that it suffices to assert that $x_j=\alpha_{j}$ for each $j$ with $\alpha_j> \alpha_{i_0}$ and $x_j=0$ for each $j$ with $\alpha_j< \alpha_{i_0}$; in view of $x_{i} \in \{0,\alpha_{i}\}$ for $i \neq i_0$, we can then simply reorder the $x_i$ for the $i$'s with $\alpha_{i}=\alpha_{i_0}$ to get a point fulfilling the claimed conditions.

First, assume there exists some $j$ with $\alpha_j> \alpha_{i_0}$ and $x_j\neq \alpha_j$, i.e., $x_j=0$.
Then, exchanging $x_j$ and $x_{i_0}$ (note $x_{i_0}\le \alpha_j$), yields a contradiction to the maximality of $\alpha_{i_0}$.

Second, assume there exists some $j$ with $\alpha_j< \alpha_{i_0}$ and $x_j\neq 0$, i.e., $x_j=\alpha_j > 0$.
By definition of $i_0$, it follows that $0< x_{i_0}< \alpha_{i_0}$. Thus, we can apply \eqref{eq_decrease}, in case $x_{i_0}< x_j$ first exchanging the two coordinates, to obtain a contradiction to the assumption that a minimum is attained in $\overline{x}$.
\end{proof}

\medskip
Note that for $G\cong C_n^r$ the bound  given by Theorem
\ref{5.4} is of the form $\daleth(H)\leq \frac 56 \mathsf D(G)+O_r(1)$. Thus, for $n$ large relative to $r$ this is  an improvement on the bound $\daleth (H) \leq \mathsf D(G)$.
\medskip
\begin{theorem} \label{5.4}
Let $H$ be a Krull monoid, $\varphi \colon H \to F = \mathcal F (P)$
a cofinal divisor homomorphism, $G = \mathcal C (\varphi)$ its class
group, and $G_P \subset G$ the set of classes containing prime
divisors. If $\exp (G) = n$ and $\mathsf r (G) = r$, then
\ber\label{sidney}\qquad \daleth (H) &\leq& \max\left\{n,\quad \frac23\mathsf
D(G_P)+\frac13\left\lfloor\left\lfloor\log_{\lfloor
n/2\rfloor+1}|G|\right\rfloor\cdot \lfloor n/2\rfloor+|G|\cdot
(\lfloor n/2\rfloor+1)^{-\lfloor\log_{\lfloor
n/2\rfloor+1}|G|\rfloor}\right\rfloor\right\}\\ &\leq& \max\left\{n,\quad \frac13\left(2\mathsf
D (G_P)+\frac{1}{2} rn+ 2^r\right)\right\}.\nn\eer
\end{theorem}

\begin{proof}
Since $\daleth (H) = \daleth (G_P)$ by Lemma \ref{3.6}.3, it suffices to show that
$\daleth (G_P)$ satisfies the given bounds. Let $U,\,V\in \mathcal
A(G_P)$ with $\max \mathsf L(UV)\geq 3$, and let
\[
z = A_1\cdot\ldots\cdot A_{r_1} B_1\cdot\ldots\cdot B_{r_2}\in
\mathsf Z(UV) \,,
\]
where $A_i,\,B_j\in \mathcal A(G_P)$ with $|A_i|\geq 3$ and
$|B_j|=2$ for all $i \in [1, r_1]$ and all $j \in [1, r_2]$, be a
factorization of $UV$ of length $|z| = \min \big( \mathsf L (UV)
\setminus \{2\} \big)$. Note $r_2\geq 2$, else $|z|\leq \frac{|UV|-2}{3}+1\leq \frac{2\mathsf D(G_P)+1}{3}$, implying \eqref{sidney} as desired (the inequality between the two bounds in Theorem \ref{5.4} will become apparent later in the proof). Our goal is to show $|z|$ is bounded above
by \eqref{sidney}. We set
\[
S=B_2 \cdot \ldots \cdot  B_{r_2}\in \mathcal B(G) \,.
\]
Observe that, for every $i \in [2, r_2]$,  $B_i$ contains one term
from $\supp(U)$ with the other from $\supp(V)$ (otherwise
$\min\{|U|,\,|V|\}=2$, contradicting $\max \mathsf L(UV) \ge 3$ in
view of Lemma \ref{5.1}). Hence we can factor $S=S_US_V$ so that
$S_U=-S_V$ with $S_U|U$ and $S_V|V$. Let $\supp (S_U) = \{g_1,
\ldots, g_{s}\}$ with the $g_i$ distinct and indexed so that $\mathsf v_{g_1}(S_U) \ge  \ldots \ge
\mathsf v_{g_s}(S_U)$. If $\mathsf v_{g_1} (S_U) \ge (n+1)/2$, then
\[
\sum_{g \in \langle g_1 \rangle} \mathsf v_g (UV) \ge \mathsf
v_{g_1} (S_U) + \mathsf v_{-g_1} (S_V) \ge n+1 \ge |\langle g_1
\rangle| + 1 \,,
\]
and Lemma \ref{5.2} implies that $|z|= \min \big( \mathsf L (UV)
\setminus \{2\} \big) \in [3, n ]$.
Therefore we may assume  $\mathsf v_{g_1}(S_U) \le \lfloor
\frac{n}{2}\rfloor$.

Suppose
\be\label{spitfireI}|S_U|>\left\lfloor\log_{\lfloor
n/2\rfloor+1}|G|\right\rfloor\cdot \lfloor n/2\rfloor+|G|\cdot
(\lfloor n/2\rfloor+1)^{-\lfloor\log_{\lfloor
n/2\rfloor+1}|G|\rfloor}-1\ee or \be|S_U|>
\frac{1}{2}nr+2^r-1.\label{spitfire2}\ee Then Lemma \ref{5.3} (applied with $\alpha=|S_U|$ and $\alpha_i=\lfloor n/2\rfloor$, and with $\alpha=|S_U|$, $\alpha_{r+1}=\max\{n/2,\,2^r-1\}$ and $\alpha_i=n/2$ for $i\neq r+1$, re-indexing the $\alpha_i$ if need be) along with $\mathsf v_{g_1}(S_U) \le \lfloor
\frac{n}{2}\rfloor\leq \frac{n}{2}$
implies that \be \label{spitfire3}\prod_{i=1}^{s} \Big( \mathsf
v_{g_i} (S_U)+1 \Big)>|G|. \ee Moreover, Lemma \ref{5.3} also
shows that the bound in \eqref{spitfireI} is at most the bound
in \eqref{spitfire2}.

Since each $g_i^{\vp_{g_i}(S_U)}$ is zero-sum free, being a
subsequence of the proper subsequence $S_U|U$, it follows that
$\{0,g_i,2g_i, \ldots, \vp_{g_i}(S_U)g_i\}$ are $\vp_{g_i}(S_U)+1$
distinct elements. Hence, in view of \eqref{spitfire3} and the
pigeonhole principle, it follows that there exists $a_i,\,b_i\in
[0,\vp_{g_i}(S_U)]$, for $i\in [1,s]$, such that, letting
$$S_A=\prod_{i=1}^s g_i^{a_i}\in \mathcal F(G_P)\quad\mbox{ and }
\quad S_B=\prod_{i=1}^s g_i^{b_i}\in \mathcal F(G_P),$$ we have
$\sigma(S_A)=\sigma(S_B)$ with $S_A\neq S_B$. Moreover, by replacing
each $a_i$ and $b_i$ with $a_i-\min\{a_i,b_i\}$ and
$b_i-\min\{a_i,b_i\}$, respectively, we may w.l.o.g. assume that
\be\label{stixy} a_i=0\quad \mbox{ or }\quad b_i=0\ee for each $i\in
[1,s]$. By their definition and in view of \eqref{stixy}, we have
\[
S_AS_B \t S_U \quad \mbox{ and }\quad (-S_B)(-S_A) \t (-S_U) =S_V
\,.
\]
From $S_A \neq S_B$, $\sigma(S_A)=\sigma(S_B)$ and $S_A|S_U$ with
$S_U$ a proper subsequence of $U\in \mathcal A(G_P)$, we conclude
that $\sigma(S_A)=\sigma(S_B)\neq 0$, and thus both $S_A$ and $S_B$
are nontrivial. Since $\sigma(S_A)=\sigma(S_B)$, we have
$\sigma(S_A(-S_B))=0$, and in view of \eqref{stixy}, the $g_i$ being
distinct and $S_A|U$ and $S_B|U$ being zero-sum free, it follows
that there is no $2$-term zero-sum subsequence in $S_A(-S_B)$. Thus,
letting $T=S_A(-S_B)$, recalling that $$S_US_V=S_U(-S_U)=S=B_2 \cdot
\ldots \cdot B_{r_2},$$ and putting all the above conclusions of
this paragraph together, we see that $T$ is a nontrivial, zero-sum
subsequence not divisible by a zero-sum sequence of length $2$ such
that $T(-T) \t B_2 \cdot \ldots \cdot B_{r_2}$. However, this leads
to factorizations $T (-T) = A_{r_1+1} \cdot \ldots \cdot A_{r'_1}$
and $S \big( (-T)T \big)^{-1} = B_2' \cdot \ldots \cdot B_{r'_2}'$,
where $A_i, B'_j \in \mathcal A (G_P)$ with $|A_i| \ge 3$ and
$|B_j'| = 2$ for all $i \in [r_1+1, r_1']$ and all $j \in [2,
r_2']$. But now the factorization
\[
z'  = A_1\cdot\ldots\cdot A_{r_1}A_{r_1+1}\cdot\ldots\cdot A_{r'_1}  B_1  B'_2\cdot\ldots\cdot
B'_{r'_2}\in \mathsf Z(UV)
\]
contradicts the minimality  of $|z| = \min \big( \mathsf L (UV)
\setminus \{2\} \big)$ (note $|z'|\geq r'_1+1\geq 3$ since $B_1|z'$ and $T$ and $-T$  were both nontrivial). So we may instead assume \be
\label{spitfire4} |S_U| \le  \left\lfloor\log_{\lfloor
n/2\rfloor+1}|G|\right\rfloor\cdot \lfloor n/2\rfloor+|G|\cdot
(\lfloor n/2\rfloor+1)^{-\lfloor\log_{\lfloor
n/2\rfloor+1}|G|\rfloor}-1  \leq  \frac{1}{2}nr+2^r-1. \ee

Now
\[
\begin{aligned}
|z| & = r_1+r_2 \le \frac{1}{3}|A_1\cdot \ldots \cdot A_{r_1}| +
\frac{1}{2} |B_1 \cdot \ldots \cdot B_{r_2}| \\ & =
\frac{1}{3}(|UV|-2|S_U|-2)+ \frac{1}{2}(2 + 2 |S_U|) \leq \frac{1}{3}
\Big( 2 \mathsf D (G_P) + |S_U| + 1 \Big) \,,
\end{aligned}
\]
which, together with \eqref{spitfire4}, implies the assertion.
\end{proof}

\medskip
As an added remark, note that the only reason to exclude the set
$B_1$ from the definition of the sequences $S$ and $S_U$ was to
ensure that $|z'|\geq 3$. However, if $r_1\geq 1$, then $|z'|\geq 3$
holds even if $B_1$ is so included. Thus the bound in \eqref{sidney}
could be improved by $-\frac{1}{3}$  in such case.

\medskip
We state one more proposition---its proof will be postponed---and
then we give the characterization of small catenary degrees.

\medskip
\begin{proposition} \label{5.5}
Let $G = C_3 \oplus C_3 \oplus C_3$. Then $\daleth (G) = \mathsf c
(G) = 4$.
\end{proposition}

\medskip
\begin{corollary} \label{5.6}
Let $H$ be a Krull monoid with class group $G$   and suppose that
every class contains a prime divisor. Then $\daleth (H)$ is finite
if and only if the catenary degree $\mathsf c (H)$ is finite if and
only if $G$ is finite. Moreover, we have
\begin{enumerate}
\item $\mathsf c (H) \leq 2$ if and only if \ $|G|\leq 2$.

\smallskip
\item $\mathsf c (H) = 3$ if and only if $G$ is isomorphic to one of the following groups\,{\rm :} $C_3, C_2 \oplus
      C_2, \mbox{ or } C_3 \oplus C_3$.

\smallskip
\item  $\mathsf c (H) = 4$ \ if and only if \
       $G$ is isomorphic to one of the following groups\,{\rm :} $ C_4, \,C_2 \oplus C_4,\, C_2 \oplus C_2 \oplus C_2,\,\mbox{ or } C_3
       \oplus C_3 \oplus C_3$.
\end{enumerate}
\end{corollary}

\begin{proof}
If $G$ is finite, then $\mathsf D (G)$ is finite (see \cite[Theorem 3.4.2]{Ge-HK06a}), and so Lemma
\ref{3.6}.3 and Theorem \ref{4.2} imply the finiteness of $\daleth (H)$ and of  $\mathsf c (H)$.
If $G$ contains elements of arbitrarily large order, then the infinitude of $\daleth(G)$ follows by Proposition \ref{4.1}.2.
And, if $G$ contains an infinite independent set, the infinitude of $\daleth(G)$ follows by Proposition \ref{4.1}.1.
In each case the infinitude of $\daleth (H)$ and $\mathsf c (H)$, thus follows by
  \eqref{crucial-inequ} and Lemma \ref{3.6}.3.
\smallskip

1. This part of the theorem is already known and included only for completeness. That $\mathsf c(H)\leq 2$ implies $|G|\leq 2$ can be found in \cite[pp. 396]{Ge-HK06a}, while $\mathsf c(H)\leq \mathsf D(G)\leq |G|$ follows from \cite[Theorem 3.4.11 and Lemmas 5.7.2 and 5.7.4]{Ge-HK06a} and implies the other direction.

\smallskip
2. See \cite[Corollary 6.4.9]{Ge-HK06a}.

\smallskip
3. Recall the comment concerning the value of $\mathsf D(G)$ after \eqref{soitis}. We may assume that $G$ is finite. Note Proposition \ref{4.1} implies $\mathsf c(G)\geq 4$ for each of the groups listed in part 3.
As noted for part 1, we have $\mathsf c(G)\leq \mathsf D(G)\leq |G|$ in general. Thus $\mathsf c(C_4)\leq 4$ and, since $\mathsf D(C_2\oplus C_2\oplus C_2)=4$, $\mathsf c(C_2\oplus C_2\oplus C_2)\leq 4$ as well.  Moreover, Corollary \ref{4.6} shows that $\mathsf c(C_2\oplus C_4)\leq \mathsf D(C_2\oplus C_4)-1=4$. Finally, $\mathsf c(C_3\oplus C_3\oplus C_3)\leq 4$ follows by Proposition \ref{5.5}. Consequently, $\mathsf c(G)=4$ for all of the groups listed in part 3.

In view of parts 1 and 2, it remains to show all other groups $G$ not listed in Corollary \ref{5.6} have $\mathsf c(G)\geq 5$. Set $\exp (G) = n$ and
$\mathsf r (G) = r$. Now Proposition \ref{4.1} shows that $\mathsf
c (G) \geq 5$ whenever $n\geq 5$ or $r\geq 4$. This leaves only $C_4\oplus C_4$, $C_4\oplus C_4\oplus C_4$, $C_2\oplus C_4\oplus C_4$ and $C_2\oplus C_2\oplus C_4$ for possible additional candidates for $\mathsf c(G)\leq 4$. However, applying Proposition \ref{4.1} to each one of these four groups shows $\mathsf c(G)\geq 5$ for each of them, completing the proof.
\end{proof}

\medskip
The remainder of this section is devoted to the proof of Proposition
\ref{5.5}, which requires some effort. Before going into details, we
would like to illustrate that geometric and combinatorial questions
in $C_3^r$ have found much attention in the literature, and our
investigations should be seen in the light of this background.  The
\ {\it Erd\H{o}s-Ginzburg-Ziv constant} \ $\mathsf s (G)$ of a
finite abelian group $G$ is the smallest integer $l \in  \N$ with
the following property:
      \begin{itemize}
      \item Every sequence $S \in \mathcal F (G)$ of length $|S| \ge l$
            has a zero-sum subsequence $T$ of length \ $|T|= \exp (G)$.
      \end{itemize}
If  $r \in \N$ and $\varphi$ is the maximal size of a cap in
AG$(r,3)$, then $\mathsf s (C_3^r) = 2 \varphi + 1$ (see
\cite[Section 5]{ E-E-G-K-R07}). The maximal size of caps in $C_3^r$
has been studied in finite geometry for decades (see
\cite{E-F-L-S02, Ed08a, Po08a}; the precise values are only known for $r
\le 6$). This shows the complexity of these combinatorial and
geometric problems. Recently, Bhowmik and Schlage-Puchta determined
the Davenport constant of $C_3 \oplus C_3 \oplus C_{3n}$. In these
investigations, they needed a detailed analysis of the group $C_3
\oplus C_3 \oplus C_3$. Building on the above results for the
Erd\H{o}s--Ginzburg--Ziv constant $\mathsf s(G)$, in particular, using that $\mathsf s
(C_3^3) = 19$, they determined the precise values of generalized
Davenport constants in $C_3^3$ (see \cite[Proposition 1]{Bh-SP07a},
and \cite{Fr-Sc10a} for more on  generalized Davenport constants).

\smallskip
We need one more definition. For an abelian group $G$ and a
sequence $S \in \mathcal F (G)$  we denote
\[
\mathsf h (S) =    \max \{ \mathsf v_g (S) \mid g \in G \} \in [0,
|S|]
   \qquad  \text{the \ {\it maximum of the multiplicities} \ of \
$S$}.
\]

We give an explicit characterization of all minimal zero-sum sequences of maximal length over $C_3^3$.
In particular, it can be seen that for this group the
Olson constant and the Strong Davenport constant do not coincide
(we do not want to go into these topics; the interested reader is
referred to Section 10 in the survey article \cite{Ga-Ge06b}).

\medskip
\begin{lemma} \label{5.7}
Let $G = C_3 \oplus C_3 \oplus C_3$ and $U \in \mathcal F (G)$. Then
the following statements are equivalent\,{\rm :}
\begin{enumerate}
  \item[(a)] $U \in \mathcal A (G)$ with $|U| = \mathsf D (G)$.

  \smallskip
  \item[(b)] There exist a basis $(e_1,e_2,e_3)$ of $G$ and $a_i, b_j \in [0,2]$ for $i \in [1,5]$ and $j \in [1,3]$ with $\sum_{i=1}^5a_i\equiv
  \sum_{j=1}^3b_j \equiv 1 \pmod{3}$
  such that
  \[
  U=e_1^2\prod_{i=1}^2(a_ie_1+e_2)\prod_{j=1}^3(a_{2+j}e_1 + b_je_2 +
  e_3) \,.
  \]
\end{enumerate}
In particular, $\mathsf h (U) = 2$ for each $U \in \mathcal A (G)$ with $|U| = \mathsf D (G)$.
\end{lemma}

\begin{proof}
Since $\mathsf D (G)=7$ (see the comments by \eqref{soitis}) it is
easily seen that statement $(b)$ implies statement $(a)$. Let $U \in
\mathcal A (G)$ with $|U| = \mathsf D (G)$. First, we assert that
$\mathsf h (U) = 2$ and, then, derive statement $(b)$ as a direct
consequence.

Since $\mathsf h(U)< \exp(G)=3$, it suffices to show $\mathsf h(U)> 1$. Assume not.
We pick some
$e_1 \in \supp (U) \subset G^{\bullet}$. Let $G = \langle e_1 \rangle \oplus K$, where $K\cong C_3\oplus C_3$ is a subgroup, and let
$\phi: G \to K$ denote the projection (with respect to this direct sum decomposition).
We set $V= e_1^{-1}U$. We observe that $\sigma(\phi(V))=0$.

We note that for each proper and nontrivial subsequence $S \t V$ with $\sigma(\phi (S))=0$, we have that $e_1\sigma(S)$ is zero-sum free, that is
\begin{equation}
\label{eq_fixedsum}
\sigma(S)=e_1.
\end{equation}
In particular, we have $\max \mathsf L (\phi(V))\le 2$ and, in combination with $\mathsf h (U)=1$, we have $0 \nmid \phi(V)$.

We assert that $\mathsf h(\phi(V))= 2$.
First, assume $\mathsf h(\phi(V))\ge 3$. This means that $V$ has a subsequence $S' = \prod_{i=1}^3(a_ie_1 + g)$ with $g \in K$ and, since $\mathsf{h}(V)=1$,  we have $\{a_1e_1,a_2e_1,a_3e_1\}= \{0,e_1,2e_1\}$ and $\sigma(S')=0$, a contradiction.
Second, assume $\mathsf h (\phi(V))=1$. Then, since $|\supp(\phi(V))|=6$ and $|K^{\bullet}|=8$,
there exist $g,h \in K$ such that $(-g)g(-h)h\t \phi (V)$, a contradiction to $\max \mathsf L (\phi(V))=2$.

So, let $g_1g_2 \t V$ with $\phi(g_1)=\phi(g_2)$, and denote this element by $e_2$. Further, let $e_3 \in K$ such that $G= \langle e_1,e_2,e_3\rangle$
and let $\phi': G \to \langle e_3 \rangle$ denote the projection (with respect to this basis).
If there exists a subsequence $T \t (g_1g_2)^{-1}V$ with $\sigma(\phi(T))= -e_2$, then
$\sigma(g_1T)$ and $\sigma(g_2T)$ are distinct elements of $\langle e_1 \rangle$, a contradiction to \eqref{eq_fixedsum}.
So, $-e_2 \notin \Sigma(\phi((g_1g_2)^{-1}V))$, which in view of $\mathsf{h} (\phi(V))<3$ and $0 \nmid \phi(V)$, implies that
$\supp(\phi((g_1g_2)^{-1}V)) \cap \langle e_2 \rangle = \emptyset$.  Since $\sigma(\phi'((g_1g_2)^{-1}V))=0$, it follows that
$\phi'((g_1g_2)^{-1}V)=e_3^2(-e_3)^2$. Let $V= g_1g_2h_1h_2f_1f_2$ such that $\phi'(h_i)=e_3$ and $\phi'(f_i)=-e_3$ for $i \in [1,2]$.
We note that $\phi(h_1 +f_1)\phi(h_2 +f_2)= 0e_2$, the only sequence of length two over $\langle e_2 \rangle$ that has sum $e_2$ yet does not have $-e_2$ as a subsum. Likewise, $\phi(h_1 +f_2)\phi(h_2 +f_1)= 0e_2$.
Thus $\phi(h_1 +f_1)=\phi(h_1 +f_2)$ or $\phi(h_1 +f_1)=\phi(h_2 +f_1)$ that is $\phi(f_1)=\phi(f_2)$ or $\phi(h_1)=\phi(h_2)$.
By symmetry, we may assume $\phi(h_1)=\phi(h_2)$. Let $j \in [1,2]$ such that $\phi(h_1 +f_j)=e_2$.
Then $\sigma( h_if_jg_1g_2) \in \langle e_1 \rangle$ for $i \in [1,2]$, yet $\sigma(h_1f_jg_1g_2)\neq \sigma(h_2f_jg_1g_2)$, as $h_1$ and $h_2$ are distinct by the assumption $\mathsf h (U)=1$. This contradicts \eqref{eq_fixedsum} and completes the argument.

It remains to obtain the more explicit characterization of $U$. Let $U=e_1^2W$ for some suitable $e_1 \in G^{\bullet}$, and let $K$ and $\phi$ as above. Similarly to \eqref{eq_fixedsum}, we see that $\phi(W)$ is a minimal zero-sum sequence over $K \cong C_3^2$. Since $\phi(W)$ has length
$5 = \mathsf{D}(C_3^2)$, it follows that $\phi(W)= e_2^2\prod_{j=1}^3(b_je_2 + e_3)$ for independent $(e_2,e_3)$ and $b_j \in [0,2]$ with $\sum_{j=1}^3b_j \equiv 1 \pmod{3}$ (cf., e.g., \cite[Example 5.8.8]{Ge-HK06a}).
Since $\sigma(W)=e_1$, the claim follows.
\end{proof}

\bigskip
\begin{proof}[Proof of Proposition \ref{5.5}]
Let $G = C_3 \oplus C_3 \oplus C_3$. Recall that $\mathsf D (G) = 7$ (see the comments by \eqref{soitis}). Thus
it  suffices to prove  $\daleth (G) \le 4$, since then combing with
Proposition \ref{4.1}.3 and Corollary \ref{4.3} yields
\[
4 \le \daleth (G) = \mathsf c (G) \le 4 \,.
\]

Suppose by contradiction that $\daleth (G)\geq 5$. Consider a counter example $U, V \in
\mathcal A (G)$ with $\max \mathsf L (UV) > 4$  and $\mathsf L (UV)\cap [3,4]=\emptyset$ such that $|U|+|V|$ is maximal.
Since $\max \mathsf L (UV)\geq 5$ and thus by Lemma \ref{5.1} $\min \{|U|,|V|\}\ge 5$, and since $\max \{|U|,|V|\}\le \mathsf D (G)=7$, we know $|U|+|V|\in [10,14]$.
Let $w=W_1\cdot \ldots \cdot W_t\in \mathsf Z(UV)$, where $t\geq 5$ and $W_i\in \mathcal A(G)$ for $i\in [1,t]$,
be a factorization of $UV$ of length at least $5$.

Note that, for some $j\in [1,t]$, say $j=1$, we must have $W_1=(-g)g$, where $g\in G^{\bullet}$, since otherwise
\[
|w| \le \lfloor \frac{|U V|}{3} \rfloor \le \lfloor \frac {14}{3}\rfloor= 4 \,,
\]
a contradiction. Since $g(-g)$ divides neither $U$ nor $V$, we may assume that $U=gU'$ and $V=(-g)V'$, where $U',\,V'\in \mathcal F(G)$ are both zero-sum free.

\medskip
\noindent

\noindent CASE 1: \,
  We have $g\notin \Sigma(U')$ or $-g\notin \Sigma(V')$,
say $g\notin \Sigma(U')$.

Then, since $-2g=g$ and $U=gU'\in \mathcal A(G)$, we have $(-g)^2U'\in \mathcal A(G)$.
Since $W_1=(-g)g$, then letting $W'_1=g^{-1}W_1(-g)^2=(-g)^3$ and
$W'_i=W_i$ for $i\in [2,t]$, we see that $w'=W'_1\cdot \ldots\cdot
W'_t\in \mathsf Z(G)$ is a factorization of $((-g)^2U')V$ with
$|w'|=t=|w|\geq 5$. As a consequence, $\max \mathsf L
((-g)^2U'V)\geq 5$, whence the maximality of $|U|+|V|$ ensures that
$((-g)^2U')V$ has a factorization
\[
z = A_1\cdot\ldots\cdot A_r\in \mathsf Z \big( (-g)^2U'V \big)
\]
with $r\in [3,4]$, where $A_i\in \mathcal A(G)$ for $i\in [1,r]$.
Note, since $-g|V$, that $\vp_{-g}((-g)^2U'V)\geq 3$.

If $(-g)^2|A_j$ for some $j\in [1,r]$, then, letting
$A'_j=A_j(-g)^{-2}g$ and $A'_i=A_i$ for $i\neq j$, gives a
factorization  $z'=A'_1\cdot\ldots\cdot A'_r\in \mathsf Z(G)$ of
$UV$ with $r\in [3,4]$ and $A'_i\in \mathcal A(G)$ for $i\in [1,r]$,
contradicting that $\mathsf L (UV)\cap [3,4]=\emptyset$. Therefore
we may assume \be\label{most1}\vp_{-g}(A_i)\leq 1\quad\mbox{ for all
} i\in [1,r].\ee As a result, since $\vp_{-g}((-g)^2U'V)\geq 3$, we
see that at least three $A_i$ contain $-g$, say w.l.o.g. $A_1$,
$A_2$ and $A_3$ with \be\label{dumphfy}|A_1|\leq |A_2|\leq |A_3|.\ee
For $i,\,j\in [1,3]$ distinct, we set
\[
B_{i,j} =(-g)^{-2}A_iA_jg\in \mathcal B(G) \,.
\]
Note that there is no $2$-term zero-sum subsequence of $B_{i,j}$
which contains $g$ as otherwise $\vp_{-g}(A_iA_j)\geq 3$, contradicting \eqref{most1}.
Consequently, \be\label{B_k-maxbound} \max \mathsf L(B_{i,j})\leq
1+\lfloor\frac{|B_{i,j}|-3}{2}\rfloor \,. \ee

\noindent CASE 1.1: \, $r=3$.

Suppose $|A_i|+|A_j|=9$ for distinct $i,\,j\in [1,3]$. Then
$|B_{i,j}|=8>\mathsf D(G)$, whence  $\min \mathsf L(B_{i,j})\geq 2$,
while \eqref{B_k-maxbound} implies $\max \mathsf L(B_{i,j})\leq 3$;
thus letting $z_B\in \mathsf Z(B_{i,j})$ be any factorization of
$B_{i,j}$, we see that $z' = z_B  A_k\in \mathsf Z(UV)$, where
$\{i,j,k\}=\{1,2,3\}$, is a factorization of $UV$ with $|z|\in
[3,4]$, contradicting $\mathsf L (UV)\cap [3,4]=\emptyset$. So we
may instead assume  \be\label{nines}|A_i|+|A_j|\neq 9\quad\mbox{ for
all distinct } i,\,j\in [1,3].\ee

Suppose $-g\in \Sigma((-g)^{-1}A_i)$ for some $i\in [1,3]$. Then,
since $\sigma((-g)^{-1}A_i)=g=-2g$, we can write $$A_i=(-g)S_1S_2$$
with $S_1,\,S_2\in \mathcal F(G)$ and $\sigma(S_1)=\sigma(S_2)=-g$.
Let $\{i,j,k\}=\{1,2,3\}$ and $\{x,y\}=\{1,2\}$. Lemma \ref{5.1}
implies $gS_x\in \mathcal A(G)$ and \be \label{stattit}\quad
(-g)^{-1}A_jS_y\in \mathcal B(G)\quad\mbox{ with }\quad \max\mathsf
L \big( (-g)^{-1}A_jS_y \big) \leq \min\{|(-g)^{-1}A_j|,\,|S_y|\}
\,.\ee Noting that $\Big((-g)^{-1}A_jS_y \Big)  \big(gS_x \big) A_k
= UV$ and letting $z_B \in \mathsf Z \big( (-g)^{-1}A_jS_y) \big)$
be any factorization of $(-g)^{-1}A_jS_y$, we see that the
factorization $z' = z_B \big( gS_x \big)  A_k \in \mathsf Z(UV)$
will contradict $\mathsf L (UV)\cap [3,4]=\emptyset$ unless
$|z_B|\geq 3$. Thus \eqref{stattit} implies $|S_y|\geq 3$ and
$|(-g)^{-1}A_j|\geq 3$. Since $y\in \{1,2\}$ and $j\in
\{1,2,3\}\setminus \{i\}$ are arbitrary, this implies first that
$|S_1|,\,|S_2|\geq 3$, whence $|A_i|\geq 7$, and second that
$|A_j|,\,|A_k|\geq 4$ for $j,\,k\neq i$. Combining these estimates,
we find that $15\leq |A_1|+|A_2|+|A_3|=|((-g)^2U')V|\leq 2\mathsf
D(G)=14,$ a contradiction. So we conclude that
\be\label{-g-notda}-g\notin \Sigma((-g)^{-1}A_i)\quad\mbox{ for all
} i\in [1,3] \,. \ee

Suppose $|A_2|\leq 4$. Let $z_B\in \mathsf Z(B_{1,3})$ be a
factorization of $B_{1,3} = \big((-g)^{-1}A_1 \big) \big(
(-g)^{-1}A_3g \big)$. In view of \eqref{-g-notda}, we see that
$(-g)^{-1}A_3g$ is zero-sum free, whence Lemma \ref{5.1} and
\eqref{dumphfy} imply $|z_B|\leq |(-g)^{-1}A_1|< |A_2|\leq 4$. Thus
$z' = z_B  A_2\in \mathsf Z(UV)$ is a factorization of $UV$ with
$|z'|\leq 4$, whence $\mathsf L (UV)\cap [3,4]=\emptyset$ implies
$|z'|=2$ and $|z_B|=1$, that is, $B_{1,3}\in \mathcal A(G)$ is an
atom. Consequently, $g^{-1}B_{1,3}=(-g)^{-2}A_1A_3$ is zero-sum
free. Hence, noting that
\[
UV  ((-g)g)^{-1} = \Big((-g)^{-2}A_1A_3 \Big)  \Big((-g)^{-1}A_2
\Big) \,.
\]
we see that Lemma \ref{5.1} implies
\[
\max \mathsf L \bigl( UV  ((-g)g)^{-1}) \bigr) \le |(-g)^{-1}A_2| <
|A_2|\leq 4 \,.
\]
which contradicts that $W  \big( (-g)g \big)^{-1}  = W_2 \cdot\ldots
\cdot W_t\in \mathsf Z \bigl( UV  ((-g)g)^{-1}) \bigr)$ is a
factorization of length $t-1=|W|-1\geq 4$. So we can instead assume
$|A_2|\geq 5$.

Observe that \be\label{no-gs}\supp((-g)^{-1}A_i)\cap \langle
g\rangle =\emptyset\quad \mbox{ for }i\in [1,3],\ee since otherwise
$\vp_{-g}(A_i)\geq 2$ or $\vp_g(UV)\geq 2$---the first contradicts
\eqref{most1}, while the the second  contradicts the supposition of
CASE 1 that $g\notin\Sigma(U')$ as $g \nmid V$. From \eqref{no-gs},
we see that $|A_1|\geq 3$, which, combined with $5\leq |A_2|\leq
|A_3|$ and $|A_1|+|A_2|+|A_3|=|((-g)^2U')V|\leq 2\mathsf D(G)=14$,
implies that
$$(|A_1|,|A_2|,|A_3|)\in \{(3,5,5),(3,5,6),(4,5,5)\}.$$ Thus, in
view of \eqref{nines}, we conclude that $|A_1|=3$ and
$|A_2|=|A_3|=5$.

Since $|B_{1,j}|=7$, for $j\in \{2,3\}$, it follows from
\eqref{B_k-maxbound} that \be\label{daisydue}B_{1,j}\in \mathcal
A(G)\quad \mbox{ for } j\in \{2,3\}\ee is an atom as otherwise $z' =
z_B  A_k\in \mathsf Z(UV)$, where $z_B\in \mathsf Z(B_{1,j})$ and
$\{1,j,k\}=\{1,2,3\}$, will contradict $\mathsf L (UV)\cap
[3,4]=\emptyset$. Since $|B_{2,3}|=9>\mathsf D(G)$, it follows from
\eqref{B_k-maxbound} that  $z' =z_B  A_1\in \mathsf Z(UV)$, for some
$z_B\in \mathsf Z(B_{2,3})$, will contradict $\mathsf L (UV)\cap
[3,4]=\emptyset$ unless all $z_B\in \mathsf Z(B_{2,3})$ have
$|z_B|=4$. Consequently, since there is no $2$-term zero-sum
containing $g$ in $B_{2,3}=(-g)^{-2}A_2A_3g$ (recall the argument
used to prove \eqref{B_k-maxbound}), we conclude that
$A_2A_3=(-g)Xa(-g)(-X) b$ for some $X=x_1x_2x_3\in \Fc(G)$ and
$a,\,b\in G$ with $$a+b=-g.$$ Thus, in view of \eqref{-g-notda}, we
find that w.l.o.g. $$A_2=(-g)Xa\quad\mbox{ and }\quad A_3=(-g)(-X)
b.$$ If $a=b$, then $2a=a+b=-g$ implies $a=g$, in contradiction to
\eqref{no-gs}. Therefore $a\neq b$.

Let $$A_1=(-g)Y\quad\mbox{ with }\quad Y=y_1y_2\in \Fc(G).$$ In view
of \eqref{daisydue}, \eqref{no-gs} and Lemma \ref{5.7}, we see that
there are terms $a'\in \supp(YXa)=\supp(B_{1,2}g^{-1})$ and $b'\in
\supp(Y(-X)b)=\supp(B_{1,3}g^{-1})$ with $$\vp_{a'}(YXa)\geq
2\quad\mbox{ and }\quad \vp_{b'}(Y(-X)b)\geq 2.$$ If $y_1=y_2$, then
$2y_1=y_1+y_2=g$ (in view of $A_1=(-g)y_1y_2$), in contradiction to
\eqref{no-gs}; if $x_i=x_j$ for $i$ and $j$ distinct, then
$x_i^2(-x_i)^2|X(-X)$, so that $x_i^2(-x_i)^2|UV$ is subsequence of
$4$ terms all from $\langle x_i\rangle$, whence Lemma \ref{5.2}
implies $UV$ has a factorization of length $3$, contradicting
$\mathsf L (UV)\cap [3,4]=\emptyset$; and if $y_i=x_j$ or $y_i=-x_j$
for some $i\in [1,2]$ and $j\in [1,3]$, then the $2$-term zero-sum
$y_i(-x_j)$ or $y_ix_i$ divides $B_{1,3}$ or $B_{1,2}$,
respectively, contradicting \eqref{daisydue}. Consequently,
$\vp_{a'}(YXa)\geq 2$ and $\vp_{b'}(Y(-X)b)\geq 2$ force $a'=a$ and
$b'=b$. Moreover, since $a\neq b$, we have $ab|XY(-X)$. Since
$a+b=-g$, we have $a^2b^2(-g)\in \mathcal B(G)$. However, noting
that there is no $2$-term zero-sum subsequence of the length $5$
zero-sum $a^2b^2(-g)$, we actually have $C = a^2b^2(-g)\in \mathcal
A(G)$. Note that $UV=g(-g)YX(-X)ab$ and $C|UV$ (in view of
$ab|XY(-X)$). Let $z_B\in \mathsf Z(UVC^{-1})$. Since
$|UVC^{-1}|=|A_1|+|A_2|+|A_3|-1-|C|=7$, we have $|z_B|\leq 3$, while
clearly $UVC^{-1}$ contains some $2$-term zero-sum subsequence from
$X(-X)$, so that $|z_B|\geq 2$. As a result, the factorization $z' =
z_B  C\in \mathsf Z(UV)$ contradicts that $\mathsf L (UV)\cap
[3,4]=\emptyset$, completing the subcase.

\noindent \noindent CASE 1.2: \, $r=4$.

If $-g\in \supp(A_4)$ as well, then we may w.l.o.g. assume
$|A_1|\leq |A_2|\leq |A_3|\leq |A_4|$, in which case
$|(-g)^2U'V|=|A_1|+|A_2|+|A_3|+|A_4|\leq 2\mathsf D(G)=14$ implies
$|B_{1,2}|\leq 5$. Thus $z'=z_B  A_3 A_4\in \mathsf Z(UV)$, where
$z_B\in \mathsf Z(B_{1,2})$, contradicts $\mathsf L (UV)\cap
[3,4]=\emptyset$ in view of \eqref{B_k-maxbound}. Therefore we may
assume $-g\notin \supp(A_4)$. Consequently, in view of \eqref{most1}
and the definition of the $A_i$, we find that $-g\notin \supp(V')$.

Since $|A_4|\geq 2$, we see that
$|(-g)^2U'V|=|A_1|+|A_2|+|A_3|+|A_4|\leq 2\mathsf D(G)=14$ implies
$|B_{1,2}|\leq 7$, with equality only possible if $|(-g)^2U'V|=14$.
However, if $|B_{1,2}|\leq 6$, then $z'=z_B  A_3 A_4\in \mathsf
Z(UV)$, where $z_B\in \mathsf Z(B_{1,2})$, contradicts $\mathsf L
(UV)\cap [3,4]=\emptyset$ in view of \eqref{B_k-maxbound}. Therefore
we indeed see that $|B_{1,2}|=7$ and $|(-g)^2U'V|=14$. As a result,
since $(-g)^2U'\in \mathcal A(G)$ implies $|U|+1=|(-g)^2U'|\leq
\mathsf D(G)=7$, and since $|V|\leq \mathsf D(G)=7$ as well, it
follows that $|V|=7$ and $|U|=6$.

Since $|V|=7= \mathsf{D} (G)$, it follows that $-g \in \Sigma(V')= G^{\bullet}$.
Thus, since $\sigma(V')=g=2(-g)$, we see that we can write $V'=S_1S_2$ with $S_1,\,S_2\in \Fc(G)$ and $\sigma(S_1)=\sigma(S_2)=-g$, and w.l.o.g. assume $|S_1|\leq |S_2|$. Then, since $|V'|=6$ and $-g\notin \supp(V')$, we infer that $2\leq |S_1|\leq 3$.

But now consider $g^{-1}U(-g)S_1\in \mathcal B(G)$ and
$((-g)S_1)^{-1}Vg\in \mathcal B(G)$. By Lemma \ref{5.1},
\[
\big(((-g)S_1)^{-1}V\big)  g \in \mathcal A(G)
\]
is an atom. Let
\[
z_B\in \mathsf Z \bigl( g^{-1}U(-g)S_1 \bigr) \,.
\]
Since $|g^{-1}U(-g)S_1|=|U|+|S_1|\geq |U|+2=8>\mathsf D(G)$, we have
$|z_B|\geq 2$. Since $g\notin \Sigma(U')=\Sigma(g^{-1}U)$ by the
supposition of CASE 1, Lemma \ref{5.1} implies $|z_B|<|(-g)S_1|\leq
4$. Thus $z'=\big( ((-g)S_1)^{-1}V \big)  z_B\in \mathsf Z(UV)$ has
$|z'|\in [3,4]$, contradicting $\mathsf L (UV)\cap [3,4]=\emptyset$
and completing CASE 1.

\medskip
\noindent CASE 2: \, We have $g\in \Sigma(U')$ and $-g\in
\Sigma(V')$.

Then, since $\sigma(U')=-g=2g$ and $\sigma(V')=g=2(-g)$,  we can write $U'=S_1S_2$ and $V'=T_1T_2$ with \
$S_1,\,S_2,\,T_1,\,T_2\in \Fc(G)$, \ $\sigma(S_1)=\sigma(S_2)=g$ \
and \ $\sigma(T_1)=\sigma(T_2)=-g$. Let $i\in \{1,2\}$ and $j\in
\{1,2\}$. Note that
\[
gT_{3-j}\in \mathcal A(G)\quad \mbox{ and }\quad (-g)S_{3-i}\in
\mathcal A(G)
\]
by Lemma \ref{5.1}. Also, $S_{i}T_{j}\in \mathcal B(G)$ and, for
$z_B\in \mathsf Z(S_{i}T_{j})$, Lemma \ref{5.1} implies
\be\label{skky}|z_B|\leq \min \{|S_{i}|,\,|T_{j}|\} \,. \ee Now
$z'=(gT_{3-j})  \bigl( (-g)S_{3-i} \bigr)  z_B \in \mathsf Z(UV)$
will contradict  $\mathsf L (UV)\cap [3,4]=\emptyset$ unless
$|z_B|\geq 3$, in which case \eqref{skky} implies $|S_{i}|\geq 3$
and $|T_{j}|\geq 3$. Since $i$ and $j$ were arbitrary, this implies
$|S_i|,\,|T_j|\geq 3$ for all $i,j\in \{1,2\}$. Hence, since
$|U|=1+|S_1|+|S_2|\leq \mathsf D(G)=7$, we see that $|S_1|=|S_2|=3$,
and likewise $|T_1|=|T_2|=3$. Thus we must have $|z_B|=3$ for all
choices of $i,j\in \{1,2\}$, which is only possible if
$S_{i}=-T_{j}$ for all choices of $i,j\in \{1,2\}$. However, this
implies $U=-V$ and, moreover, that $\vp_{-x}(T_{i})\geq 1$ for $i\in
[1,2]$ and $x\in \supp(S_1S_2)$. Consequently, letting $x\in
\supp(S_1S_2)$, we see that $\vp_{-x}(V)\geq 2$, whence $U=-V$
implies $\vp_x(U)\geq 2$. Thus $x^2(-x)^2 \t UV$ is a subsequence of
$4$ terms all from $\langle x\rangle$, whence Lemma \ref{5.2}
implies $UV$ has a factorization of length $3$, contradicting
$\mathsf L (UV)\cap [3,4]=\emptyset$ and completing CASE 2 and the
proof.
\end{proof}

\providecommand{\bysame}{\leavevmode\hbox
to3em{\hrulefill}\thinspace}
\providecommand{\MR}{\relax\ifhmode\unskip\space\fi MR }
\providecommand{\MRhref}[2]{%
  \href{http://www.ams.org/mathscinet-getitem?mr=#1}{#2}
} \providecommand{\href}[2]{#2}

\end{document}